\documentclass{amsart}

\usepackage{mathrsfs}

\ifx \relax
      \else
      \fi

\usepackage{amssymb}
\usepackage{amsthm}

\input setup.sty
\input macros220810.sty

\usepackage[hidelinks]{hyperref}

\begin{document}
\makeatletter\def\shfiuwefootnote{\gdef\@thefnmark{}\@footnotetext}\makeatother\shfiuwefootnote{Version 2022-10-14\_2. See \url{https://shelah.logic.at/papers/363/} for possible updates.}



\title{On Spectrum of $\kappa$--Resplendent Models}
\author{Saharon Shelah}
\thanks{This was supposed to be Ch. V of the book ``Non-structure" 
and probably will be if it materializes. Was circulated around 1990.
First typed in 1988. The author would like to thank ISF-BSF for partially supporting this research by grant with Maryanthe number NSF 2051825, BSF 3013005232. \\
References like [Sh:950, Th0.2=Ly5] mean that the internal label of Th0.2 is y5 in Sh:950.
The reader should note that the version in my website is usually more up-to-date than the one in arXiv.
This is publication number    
363
in Saharon Shelah's list.\\
On the old versions, the author expresses gratitude for the partial support of the Binational Science Foundation in this research and thanks Alice Leonhardt for her careful and beautiful typing. 
In new versions, the author thanks an individual who wishes to remain anonymous for funding typing services, and thanks Matt Grimes for the careful and beautiful typing.
}


\makeatletter
\@namedef{subjclassname@2020}{\textup{2020} Mathematics Subject Classification}
\makeatother
\subjclass[2020]{FILL}
\date{2022-08-15}

\begin{abstract}
We prove that some natural ``outside" property 
of counting models up to isomorphism is equivalent 
(for a first order class) to being stable.

For a model, being resplendent is a strengthening of being $ \kappa $-saturated.
Restricting ourselves to the case $\kappa > |T|$ for transparency,
to say a model $ M $ is $ \kappa $-resplendent means:
\begin{quote}
when we expand $ M $ by $ < \kappa $ individual constants $ \LL c_i : i < \alpha \RR $,
 if $(M, c_i)_{ < \alpha }$ has an elementary extension expandable to be a model of $ T'$
where ${\rm Th}((M, c_i)_{i < \alpha} ) \subseteq T'$, $|T'| < \kappa$
then already $(M, c_i)_{i < \alpha}$ can be expanded to a model of $ T'$ .
\end{quote}
Trivially,  any saturated model of cardinality $\lambda$ 
is $\lambda$-resplendent. 
We ask: how many $\kappa $-resplendent models
of a (first order complete) theory $ T $ 
of cardinality $ \lambda $ are there?
We restrict ourselves to cardinals $\lambda = \lambda^\kappa + 2^{|T|}$
and ignore the case $\lambda = \lambda^{<\kappa} + |T| < \lambda^\kappa$.
Then we get a complete and satisfying answer: this depends
only on $ T $ being stable or unstable.
In this case proving that for stable $ T $
we get few, is not hard; in fact, every resplendent model of 
$ T $ is saturated hence it is determined by its
cardinality up to isomorphism. The inverse is
more problematic because naturally we have to use
Skolem functions with  any $ \alpha < \kappa $ places.
Normally we use relevant partition theorems
(Ramsey theorem or Erd\H{o}s-Rado theorem),  
but in our case the relevant partitions theorems fail
so we have to be careful.
\end{abstract}

\maketitle


\setcounter{section}{-1}


\section{ Introduction}
Our main conclusion speaks on stability of first order theories, but the
major (and the interesting) part of the proof has little to do with it and
can be read without knowledge of classification theory (only the short proof
of \ref{1.6} uses it), except the meaning of $\kappa<\kappa(T)$ which
we can take as the property we use, see 
inside
\ref{1.17}(1)  here (or  see   
     \cite[1.5(2)]{Sh:E59}  
or \cite{Sh:c}). 
 The point is to construct a model in which, for some infinite
sequences of elements, we have appropriate automorphisms, 
so we need to use ``Skolem" functions with infinitely many places.
Now having functions with infinite arity makes obtaining
models generated by indiscernibles harder. More specifically, the theory of
the Skolemizing functions witnessing resplendence
for $(M,\bar{b})$ is not continuous in ${\rm Th}(M,\bar{b})$. Hence we use a 
weaker version of indiscernibility though having a linear order is usually
a very strong requirement (see \cite[\S3]{Sh:E59}), 
in our proof we use it as if we only have trees (with $\kappa$ levels).

In \cite{Sh:a} or \cite[VI 5.3--5.6]{Sh:c} we characterized first order
$T$ and cardinals $\lambda$ such that for some first order complete
$T_1$, $T\subseteq T_1$, $|T_1|=\lambda$ and any
$\tau (T)$--reduct of a model of $T_1$ is saturated.

In \cite{Sh:225} we find the spectrum of strongly
$\aleph_\eps$--saturated models, but have nothing comparable for
strongly $\aleph_1$--saturated ones (on better computation of the numbers
see \cite{Sh:225a}, and more in 
     \cite[3.2]{Sh:331}). Our interest
was:
\begin{enumerate}
\item[(A)] an instance of complete classification for an ``outside''
question: 
the question here is the function giving  the number, up to isomorphism,	 
of $\kappa $-resplendent models of a (first order complete) theory $ T $ 
as a function of the cardinality; we concentrate on the case 
$ \lambda = \lambda^\kappa + 2^{|T|}$. 
\item[(B)] an ``external'' definition of stability which happens to be the
dividing line.
\end{enumerate}

Earlier we had such an equivalent 
``external" definition of stability 
by saturation of ultra-powers, i.e. Keisler order, 
see \cite{Sh:c}.  
Baldwin had told me he was writing a paper on resplendent
models: for
$\aleph_0$--stable one there are few ($\leq 2^{\aleph_0}$) such models in
any cardinality; and for $\tT$ not superstable --- there are $2^\lambda$
models of cardinality $\lambda$ (up to isomorphism).
\medskip

Note that resplendent
models are strongly $\aleph_0$--homogeneous and really
the non-structure are related.
The reader may thank Rami Grossberg for urging the author 
to add more explanation to \ref{1.7}.

\begin{Notation}\label{z2}
1) For a model $M$ and $\bar c \in {}^\alpha\! M$, let $(M,\bar c)$ 
(or $(M, c_i)_{i < \alpha}$) be $M$ expanded by the individual constants 
$c_i$ for $i < \alpha$.

2) $\bar x_{[u]} = \LL x_\eps : \eps \in u \RR$
\end{Notation}

\begin{Definition}\label{z14}
For a complete first order theory we defined $\kappa(T)$, an infinite cardinal 
or $\infty$, minimal such that:
\begin{itemize}
    \item $\kappa < \kappa(T)$ iff some $\bar \varphi$ witnesses it, which means:
    \begin{enumerate}
        \item $\bar \varphi = \LL \varphi_\alpha(x, \bar y_\alpha) : \alpha < \kappa \RR$
        
        \item for any $\lambda$, for some model $M$ of $\tT$ and sequence 
        $\LL \bar a_\eta : \eta\in {}^{\kappa\geq}\lambda\RR$ with 
        $\bar a_\eta \in {}^{\lh(\bar y_\alpha)}\!M$ we have: 
        if $\eps<\kappa$, $\eta\in {}^\kappa\!\lambda$, $\alpha<\lambda$ then
        $$M \models \varphi_\eps[a_\eta, \bar a_{\eta\rest\eps\conc\LL\alpha\RR}]^{{\rm if}(\alpha=\eta(\eps))}.$$
    \end{enumerate}
\end{itemize}
\end{Definition}

\begin{Claim}\label{z17}
In \ref{z14}, for our purposes, we can demand $\bar y_\alpha = \LL y \RR$.
\end{Claim}

\begin{PROOF}{\ref{z17}}
Letting $\gC$ be a model of $T$, use $\Th(\gC^\cf)$ (see \cite[Ch.III]{Sh:c}).

If $\cf(\kappa) > \aleph_0$, without loss of generality $\LL \lh(\bar y_\alpha) : \alpha < \kappa \RR$ 
is constant (call this value $n$) so we can use essentially ${}^n\gC$. This holds for the main 
lemma \ref{1.7} except when $\kappa = \aleph_0$, in which case the results are easy.
\end{PROOF}

\begin{Definition}\label{z20}
The sequences $\bar \eta$, $\bar \nu$ from ${}^{\kappa \geq}\!\mu$ are called \emph{similar} when:
\begin{enumerate}
    \item $\lh(\bar \eta) = \lh(\bar \nu)$ (call this $n$)
    
    \item if $\ell < n$ then $\lh(\eta_\ell) = \lh(\nu_\ell)$
    
    \item if $k,\ell < n$ \underline{then}
    \begin{enumerate}
        \item $\eta_k \lhd \eta_\ell$ iff $\nu_k \lhd \nu_\ell$
    
        \item $\lh(\eta_k \cap \eta_\ell) = \lh(\nu_k \cap \nu_\ell)$
    
        \item if $\eps = \lh(\eta_k \cap \eta_\ell)$ is $<$ both 
        $\lh(\eta_k)$ and $\lh(\eta_\ell)$ then 
        $$ \eta_k(\eps) < \eta_\ell(\eps) \Leftrightarrow \nu_k(\eps) < \nu_\ell(\eps)$$

    \end{enumerate}
\end{enumerate}
\end{Definition}

\section{Resplendency}
Our aim is to prove \ref{1.2} below (``$\kappa$--resplendent'' is defined in
\ref{1.3}).

\begin{Convention}
\label{1.1}
{\rm
$\tT$ is a fixed first order complete theory; recall that 
$\tau(T) = \tau_T $, 
$\tau(M) = \tau _ M $ 
is the vocabulary of $\tT$, $M$ respectively and
$\bbL$ is first order logic, 
so $\bbL_\tau\equiv \bbL (\tau)$ is
the first order language with vocabulary $\tau$.
}
\end{Convention}

We show here

\begin{Theorem}
\label{1.2}
The following are equivalent (see Definition \ref{1.3} below) for a regular
uncountable $\kappa$:
\begin{enumerate}
\item[(i)]   $\kappa<\kappa(\tT)$, see e.g. \ref{z14} or~\ref{1.17}(1),
\item[(ii)]  there is a non-saturated $\kappa$--resplendent
model of $\tT$
(see Definition \ref{1.3} below),
\item[(iii)] for every $\lambda=\lambda^\kappa\geq 2^{|\tT|}$, $\tT$
has $>\lambda$ non-isomorphic $\kappa$--resplendent models,
\item[(iv)]  for every $\lambda=\lambda^\kappa\geq 2^{|\tT|}$, $\tT$
has $2^\lambda$ non-isomorphic $\kappa$--resplendent
models.
\end{enumerate}
\end{Theorem}

\begin{PROOF}{\ref{1.2}}
The implication ${\bf (i)}\ \Rightarrow\ {\bf (iii)}$ follows from 
the main Lemma \ref{1.7} below; the implication  ${\bf (iii)}\ \Rightarrow\ {\bf (ii)}$ 
is trivial (as any two saturated models of $T$ of the same cardinality are isomorphic), 
and ${\bf (ii)}\ \Rightarrow\ {\bf (i)}$ follows from \ref{1.6} below. Trivially, 
${\bf (iv)}\ \Rightarrow\ {\bf (iii)}$, and lastly ${\bf (i)}\ \Rightarrow\ {\bf (iv)}$ 
by \ref{3.1}+\ref{1.33}.
\end{PROOF} 

\begin{Remark}
\label{1.2A}
{\rm
\begin{enumerate}
    \item If we omit condition (iv) we save \S3 as well as the dependency on a theorem 
    from \cite{Sh:309} 
    using only an easy relative. This is sufficient for having an outside 
    characterization of $T$ be stable.
    
    \item In the proof the main point is ${\bf (i)}\ \Rightarrow\ {\bf (iii)}$
    (and ${\bf (i)}\ \Rightarrow\ {\bf (iv)}$, i.e., the non-structure part).
    
    \item Remember: $\tT$ is unstable iff $\kappa(\tT)=\infty$.
    
    \item Notice that every saturated model $M$ is $\|M\|$--resplendent
    (see \ref{1.3}(2) below). Actually a little more.
\end{enumerate}
}
\end{Remark}

\begin{Definition}
\label{1.3}
\begin{enumerate}
\item[(1)] A model $M$ is $(\kappa,\ell)$--resplendent
(where $\ell=0,1,2,3$)
if: \underline{for every} elementary extension $N$ of $M$ and expansion
$N_1$ of $N$ satisfying $|\tau(N_1)\setminus\tau(N)|<\kappa$ and $\alpha<
\kappa$, $c_i\in M$ for $i<\alpha$ and $\tT_1\subseteq {\rm Th}(N_1,c_i)_{i<
\alpha}$ satisfying $(\ast)^\ell_{\tT_1}$ below, \underline{there is} an
expansion $\left(M_1,c_i\right)_{i<\alpha}$ of $\left(M,c_i\right)_{i<
\alpha}$ to a model of $\tT_1$,
when: 
\begin{enumerate}
    \item[$(\ast)^\ell_{\tT_1}$]
    \begin{description}
        \item[\underline{Case 0}] $\ell=0$:\quad $|\tT_1|<\kappa$,
        \item[\underline{Case 1}] $\ell=1$:\quad for some 
        $\tau'\subseteq \tau(N_1)$, $|\tau'|<\kappa$ and 
        $$\tT_1\subseteq\bbL (\tau'\cup \{c_i:i<\alpha\})$$
        
        \item[\underline{Case 2}] $\ell=2$:\quad $\tT_1$ is
        $\kappa$--recursive (see \ref{1.3}(4) below),
        
        \item[\underline{Case 3}] $\ell=3$:\quad 
        $\tT_1={\rm Th}\left(N_1, c_i\right)_{i<\alpha}$ (but remember that 
        $N_1$ has only $<\kappa$ relations and functions not of $M$).
    \end{description}
\end{enumerate}
\item[(2)] $\kappa $-resplendent
means $(\kappa ,3)$-resplendent.
\item[(3)] Assume $M$ is a model of $\tT$, $\bar c\in {}^{\kappa>}|M|$ and
$M_{\bar c}$ is an expansion of $(M,\bar c)$.
We say that $M_{\bar c}$
witnesses $(\kappa,\ell)$--resplendence for $\bar c$ in $M$,
\underline{when}:

for every first order $\tT_{1}$ such that
\[{\rm Th}(M,\bar c)\subseteq \tT_1\quad \&\quad |\tau(\tT_1)
\setminus\tau(\tT)|<\kappa\]
and $(\ast)^\ell_{\tT_1}$ holds, we have:

$M_{\bar c}$ is a model of $\tT_1$ up to renaming the symbols in
$\tau(\tT_{1})\setminus \tau(M,\bar c)$.
\item[(4)] For $M,N_1$, $\LL c_i:i<\alpha\RR$ and $\tT_1\subseteq
{\rm Th}(N_1,c_i)_{i<\alpha}$ as in part (1), $\tT_1$ is
$\kappa$--recursive \underline{when}:
\begin{enumerate}
    \item[(a)] $\kappa=\aleph_0$ and $\tT_1$ is recursive (assuming the 
    vocabulary of $\tT$ is represented in a recursive way) \underline{or}
    
    \item[(b)] $\kappa>\aleph_0$ and for some $\tau^\ast\subseteq\tau(N_1)$, 
    $|\tau^\ast|<\kappa$ the following holds:

    \underline{if}\ $\varphi_\ell(x_0,\ldots,x_{n-1})\in \bbL(\tau')$ for
    $\ell=1,2$ and there is an automorphism $\pi$ of $\tau'$ (see parts (9), (10)),
    where $\tau^\ast\subseteq\tau'\subseteq\tau(N_1)$ such that
    $\pi$ is the identity on $\tau^\ast$ and $\hat{\pi}(\varphi_1)=\varphi_2$ and
    $\beta_0<\beta_1<\ldots<\alpha$ \underline{then}
\[
\varphi_1(c_{\beta_0},c_{\beta_1},\ldots)\in \tT_1\quad\mbox{ iff }
\quad\varphi_2(c_{\beta_0},c_{\beta_1},\ldots)\in \tT_1.
\]
\end{enumerate}
\item[(5)] We say $f$ is an $(M,N)$-elementary mapping \underline{when} $f$ is a
partial one-to-one function from $M$ to $N$, $\tau(M)=\tau(N)$ and for every
$\varphi(x_0,\dots,x_{n-1})\in \bbL (\tau_M)$ and $a_0,\ldots,a_{n-1}\in M$ we have:
\[
M\models\varphi[a_0,\ldots,a_{n-1}]\quad\mbox{ iff }\quad N\models
\varphi[f(a_0),\ldots,f(a_{n-1})].
\]

\item[(6)] $f$ is an $M$--elementary mapping \underline{if} it is an $(M,M)$--elementary
mapping.
\item[(7)] $M$ is $\kappa$--homogeneous \underline{if}: 

\underline{for any} $M$--elementary mapping $f$ with $|\Dom(f)|<\kappa$ and
$a\in M$ \underline{there is} an $M$--elementary mapping $g$ such that:
\[f\subseteq g,\qquad \Dom(g)=\{a\}\cup\Dom(f).\]
\item[(8)] $M$ is strongly $\kappa$--homogeneous 
\underline{if}\  
for any $M$-elementary mapping $f$ with $|\Dom(f)|<\kappa$ there is an 
automorphism $g$ of $M$ such that $f\subseteq g$.

\item[(9)] Let $\tau_1\subseteq \tau_2$ be vocabularies. We say that $\pi$ is an 
automorphism of $\tau_2$ over $\tau_1$ \underline{when}: $\pi$ is a permutation of $\tau_2$,
$\pi$ maps any predicate $P\in \tau_2$ to a predicate of $\tau_2$
with the same arity, $\pi$ maps any function symbol
$F \in \tau_2$ to a function symbol of $\tau_2$ of the same
arity and $\pi\rest \tau_1$ is the identity.

\item[(10)] For $\pi, \tau_2$ as in part (9) let $\hat{\pi}$ be the permutation 
of the set of formulas in the vocabulary $\tau_2$ which $\pi$ induces.
\end{enumerate}
\end{Definition}

Note:

\begin{Fact}
\label{1.4A}
\begin{enumerate}
\item If $\tau=\tau(M)$, and
$$\Big[\tau'\subseteq\tau\ \&\ |\tau'|<\kappa\ \Rightarrow\  M
\rest\tau'\mbox{ is saturated }\Big]
$$
\underline{then} $M$ is $(\kappa,1)$--resplendent.
\item If $M$ is saturated of cardinality $\lambda$ \underline{then} $M$ is
$\lambda$-resplendent.
\end{enumerate}
\end{Fact}

\begin{PROOF}{\ref{1.4A}}
Easy, e.g., see \cite{Sh:a} and not used here elsewhere.
\end{PROOF}

\begin{Example}
\label{1.4}
{\rm
There is, for each regular $\kappa$, a theory $\tT_\kappa$ such that:
\begin{enumerate}
    \item[(A)] $\tT_\kappa$ is superstable of cardinality $\kappa$,
    
    \item[(B)] for $\lambda\geq\kappa$, $\tT_\kappa$ has $2^\lambda$ 
    non-isomorphic $(\kappa,1)$--resplendent models.
\end{enumerate}
}
\end{Example}

\par \noindent

\begin{PROOF}{\ref{1.4}}
Let $A_0=\{\kappa\setminus (i+1): i<\kappa\}$ and $A_1= A_0
\cup \{\varnothing \}$. For every linear order $I$ of cardinality $\lambda\geq
\kappa$ we define a model $M_I$:

\noindent its universe is
$$
I\cup\Big\lbrace\LL s,t,i,x\RR : s\in I,\ t\in I,\ i<\lambda,\ x\in A_1
\mbox{ and }\big[ I \models s < t\ \Rightarrow\ x\in A_0\big]\Big\rbrace,
$$
(And of course, without loss of generality, no quadruple $\LL s,t,i,x\RR$ 
as above belongs to $I$.) Its relations are:
\[\begin{array}{lcl}
P&=&I,\\
R&=&\left\lbrace\big\LL s,t,\LL s,t,i,x\RR\big\RR:s\in I,\
t\in I,\ \LL s,t,i,x\RR\in |M_I|\setminus P\right\rbrace,\\
Q_\alpha&=&\left\lbrace\LL s,t,i,x\RR:\LL 
s,t,i,x\RR\in  
|M_I |  
\setminus P,\ \alpha\in x\right\rbrace\quad \mbox{ for }\alpha<\kappa.
  \end{array}\]
In order to have the elimination of quantifiers we also have two unary
functions $F_1$, $F_2$ defined by:
\[\begin{array}{rcl}
s\in I&\Rightarrow & F_1(s)=F_2(s)=s,\\
\LL s,t,i,x\RR\in  
|M_I |  \setminus  
I&\Rightarrow& F_1(\LL s,t,i,x
\RR)=s\ \&\ F_2(\LL s,t,i,x\RR)=t.
  \end{array}\]
It is easy to see that:
\begin{enumerate}
    \item[(a)] In $M_I$, the formula
\[
P(x)\ \&\ P(y)\ \&\ (\exists z)(R(x,y,z)\ \&\ \bigwedge\limits_{\alpha < \kappa}\neg Q_\alpha(z))
\]
    linearly orders $P^{M_I}$; in fact, it defines $<_I$.
    
    \item[(b)] ${\rm Th}(M_I)$ has elimination of quantifiers.

    \item[(c)] If $\tau\subseteq\tau(M_I)$, $|\tau|<\kappa$ then 
    $M_I\rest\tau$ is saturated.
    \item[(d)] ${\rm Th}(M_I)$ does not depend on $I$ (as long as it is
    infinite) and we call it $\tT_\kappa$.
    
    \item[(e)] $\tT_\kappa$ is superstable.
\end{enumerate}
Hence: $\tT_\kappa={\rm Th}(M_I)$ is superstable, does not depend on $I$,
and
\[M_I\cong M_J\quad \mbox{\underline{iff}}\quad I\cong J,\]
and by \ref{1.4A} $M_I$ is $(\kappa,1)$--resplendent.

This suffices for part (A) of the claim. By clause (e) above, part (B) of the claim 
will follow by \cite[IV,\S3]{Sh:300} (or better, \cite[\S3]{Sh:E59}).
\end{PROOF} 

\begin{Fact}
\label{1.5}
\begin{enumerate}[(1)]
    \item for $\ell = 1,2$, $M$ is $(\kappa,3)$--resplendent implies $M$ is 
    $(\kappa,\ell)$--resplendent, which implies $M$ is 
    $(\kappa,0)$--resplendent.
    
    \item $M$ is $(\kappa,0)$--resplendent implies $M$ is $\kappa$-compact.
    
    \item $M$ is $(\kappa,2)$--resplendent implies $M$ is $\kappa$-homogeneous,
    even strongly $\kappa$-homogeneous (see Definition \ref{1.3}(7),(8)).
    
    \item If $M$ is $(\kappa,2)$--resplendent $\kappa>\aleph_0$ and 
    $\{\bar a_n : n < \omega\}$ is an indiscernible set in $|M|$,  
    \underline{then} it can be extended to an indiscernible set of 
    cardinality $\|M\|$ (similarly for sequences).
    
    \item $M$ is $(\kappa,3)$--resplendent implies $M$ is $\kappa$-saturated.
    
    \item If $\kappa>|\tT|$ \underline{then} the notions of \ref{1.3}
    ``$(\kappa,\ell)$--resplendent'' for $\ell= 0,1,2,3$, are equivalent.
\end{enumerate}
\end{Fact}

\begin{PROOF}{\ref{1.5}}
Straightforward: for example,\\
(3)\quad For given $a_i,b_i\in M$ (for $i<\alpha$, where $\alpha<\kappa$) use
\[\begin{array}{ll}
\tT_1=&\{G(a_i) = b_i : i < \alpha\}\ \cup\\
&\{(\forall x,y)(G(x) = G(y)\ \Rightarrow\ x = y),\ (\forall x)(\exists y)(G(y) = x)\}\ \cup\\
& \{(\forall x_0,\ldots,x_{n-1})[R(x_0,\ldots,x_{n-1})\equiv R(G(x_0),
\ldots,G(x_{n-1}))]:\\
&\qquad\qquad R\mbox{ an $n$-place predicate of }\tau(M)\}\ \cup\\
&\{(\forall x_0,\ldots,x_{n-1})[F(G(x_0),\ldots)=G(F(x_0,\ldots))]:\\
&\qquad\qquad F\mbox{ an $n$-place function symbol of }\tau(M)\}.
\end{array}\]
(4) For notational simplicity let $\bar a_n=a_n$. Let $\tT_1$ be, with
$P$ a unary predicate, $g$ a unary function symbol,
\[\begin{array}{l}
\{\mbox{ ``$G$ is a one-to-one function into $P$'' }\}\cup\{P(a_n):n<
\omega\}\ \cup\ \\
\begin{array}{lr}
\big\{(\forall x_0,\ldots,x_{n-1})&\Big[\bigwedge\limits_{\ell<n}
P(x_\ell)\ \&\ \bigwedge\limits_{\ell<m<n} x_\ell\neq x_m\ \&\
\varphi[a_0,\ldots,a_{n-1}]\qquad\\
&\Rightarrow\ \varphi(x_0,\ldots,x_{n-1})\Big]:\quad\\
&\varphi(x_0,\ldots,x_{n-1})\in \bbL(\tau(M))\big\}
\end{array}
\end{array}\]
\end{PROOF} 

\begin{Conclusion}
\label{1.6}
If $M$ is $\kappa$--resplendent, $\kappa\geq\kappa(\tT)+\aleph_1$
\underline{then} $M$ is saturated.
\end{Conclusion}

\begin{PROOF}{\ref{1.6}}
By \ref{1.5}(5) $M$ is $\kappa$--saturated, so without loss of
generality $\|M\| \geq \kappa$ (and $\kappa \geq \aleph_1$). Hence, by \cite{Sh:a} or
\cite[III,3.10(1),p.107]{Sh:c}, it is enough to prove: for $\seqI$
an infinite indiscernible $\subseteq M$, $\dim(\seqI,M)=\|M\|$. But this
follows by \ref{1.5}(4). 
\end{PROOF} 

\begin{Main Lemma}
\label{1.7}
Suppose that $\kappa=\cf(\kappa)<\kappa(\tT)$ (for example, $\tT$
unstable, $\kappa$ regular) and $\lambda=\lambda^\kappa+2^{|\tT|}$.
\underline{then} $\tT$ has $>\lambda$ pairwise non-isomorphic
$\kappa$--resplendent models of cardinality $\lambda$.
\end{Main Lemma}

Before embarking on the proof, we give some explanations.

\begin{Remark}
\label{1.7A}
{\rm
\begin{enumerate}[(1)]
    \item We conjecture that we can weaken in \ref{1.7} the hypothesis
    ``$\lambda=\lambda^\kappa+2^{|\tT|}$'' to 
    ``$\lambda = \lambda^{<\kappa} + 2^{|\tT|}$''. This holds for many 
    $\lambda$-s (see \cite[\S2]{Sh:309}), and probably for all, but 
    we have not looked at this. See \S3.
    
    \item We naturally try to imitate \cite{Sh:a}, \cite[VII,\S2,
    VIII,\S2]{Sh:c} or \cite[\S3]{Sh:E59},\cite{Sh:331}.   
    In the proof of the theorem, the difficulty is that while expanding 
    to take care of resplendence, we naturally will use Skolem functions 
    with infinite arity, and so we cannot use compactness so easily.

    If the indiscernibility is not clear, the reader may look again at
    \cite[VII,\S2]{Sh:a} or \cite[VII,{\S}2]{Sh:c}, (tree indiscernibility). 
    We get below first a weaker version of indiscernibility,
    as it is simpler to get it, and is totally harmless if we would like 
    just to get $>\lambda$ non-isomorphic models by the old version 
    \cite[III,4.2(2)]{Sh:300} or the new \cite[\S2]{Sh:309}  
\end{enumerate}
}
\end{Remark}

\begin{Explanation}
\label{1.7B}
{\rm
Note that the problem is having to deal with sequences of $<\kappa$ elements
$\bar b = \LL b_i : i < \eps \RR$, $\eps$ infinite. The need to
deal with such $\bar b$ with \underline{all theories} of small vocabulary is
not serious -- there is a ``universal one'' though possibly of larger
cardinality, i.e., if $M\models \tT$, $b_i\in M$ for $i < \eps$,
$\eps < \kappa$, we can find a f.o.~theory $\tT_2 = T_2 (\bar b)$ satisfying 
${\rm Th}(M,b_i)_{i<\eps}\subseteq \tT_2$, 
$|\tT_2|\leq (2^{|\tT|+ |\eps|})^{<\kappa}$ such that:
\begin{itemize}
    \item if ${\rm Th}(M,b_i)_{i<\eps}\subseteq \tT'$ and
    $|\tau(\tT')\setminus\tau(\tT)\setminus\{b_i:i<\eps\}|<\kappa$
    \underline{then} renaming the predicates and function symbols outside $T$,
    we get $\tT'\subseteq \tT_2(\bar b)$.
\end{itemize}
This is possible by the Robinson consistency lemma. Let us give more details.
}
\end{Explanation}

\begin{Claim}
\label{1.7C}
\begin{enumerate}
    \item[(1)] Let $M_0$ be a model, $\tau_0 = \tau(M_0)$, $\bar b = \LL b_i : i < \eps \RR$ 
    where $b_i\in M_0$ for $i<\eps$ and $\theta\geq\aleph_0$ be a cardinal. 
    Let $\tau_1=\tau_0\cup\{b_i:i<\eps\}$ so $M_1=(M_0,b_i)_{i<\eps}$ is a
    $\tau_1$--model. \underline{then} there is a theory $\tT_2=\tT_2[\bar b]=\tT_2[\bar{b},M]$, 
    depending \underline{only} on $\tau_0$, $\tau_1$ and ${\rm Th}(M_1)$, i.e., 
    essentially on $\tp(\LL b_i : i < \eps \RR, \varnothing, M_0)$ such that:
    \begin{enumerate}
        \item[(a)] $\tau_2=\tau(\tT_2)=\tau(\eps,\tau_0)$ extends $\tau_1$ and
        has cardinality $\leq 2^{|\tau_1|+\theta+|\eps|}$,
        
        \item[(b)] for every $M_2$, $\tT'$, the model
        $M_2$ is expandable to a model of $\tT'$,
        \underline{when}:\\
        \begin{enumerate}
            \item[$(\alpha)$] $M_2$ is a $\tau_1$--model,
            
            \item[$(\beta)$]  $M_2$ can be expanded to a model of $\tT_2$,

            \item[$(\gamma)$] ${\rm Th}(M_2)\subseteq \tT'$, equivalently some 
            elementary extension of $M_2$ is expandable to a model of $\tT'$,

            \item[$(\delta)$] $\tT'$ is f.o.~and $|\tau(\tT')\setminus\tau(M_2)| \leq\theta$,
        \end{enumerate}
        
        \item[(a)$^+$] if $\theta>|\tT|+|\eps|$ then $|\tau_2|\leq 2^{<\theta}$ is enough.
    \end{enumerate}
    
    \item[(2)] If in part (1), sub-clause $(\delta)$ of clause (b) is weakened to:
    \begin{enumerate}
        \item[$(\delta)_2$] $\tT'$ is f.o., and $|\tau(\tT')\setminus\tau(M_2)|<\theta$,
    \end{enumerate}
    \underline{then} we can strengthen (a) to
    \begin{enumerate}
        \item[(a)$_2$] $\tau_2=\tau(\tT_2)$ extends $\tau_1$ and has cardinality
        $\leq\sum\limits_{\mu<\theta}2^{|\tau_1|+\mu+\aleph_0+|\eps|}$,

        \item[(a)$^+_2$] if $\theta>(|\tT|+|\eps|)^+$ then 
        $|\tau_2|\leq\sum\limits_{\mu<\theta}2^{<\mu}$ is enough.
    \end{enumerate}
\end{enumerate}
\end{Claim}

\begin{PROOF}{\ref{1.7C}}
 1) We ignore function symbols and individual constants as we
can replace them by predicates. Let
\[
\begin{array}{lr}
 \bbT = \big\{\tT' : &\tT' \mbox{ f.o.~complete theory, }{\rm Th}(M_1)
 \subseteq \tT'\mbox{ and\quad }\\
 &\tau(\tT')\setminus\tau(M_1)\mbox{ has cardinality }\leq \theta \big\}.
\end{array}
\]
This is a class; we say that $\tT'$, $\tT''\in\bbT$ are
isomorphic over ${\rm Th}(M_1)$ (see \cite{Sh:8}) 
\underline{when} 
there is a
function $\bfh$ satisfying:
\begin{enumerate}
    \item[(a)] $\bfh$ is one-to-one,
    
    \item[(b)] $\Dom(\bfh)=\tau(\tT')$,
    
    \item[(c)] $\Rang(\bfh)=\tau(\tT'')$,
    
    \item[(d)] $\bfh$ preserves arity (i.e., the number of places, and of course
    being predicate/function symbols),
    
    \item[(e)] $\bfh\rest(\tau(M_1))= {\rm identity}$,
    
    \item[(f)] for a f.o.~sentence $\psi=\psi(R_1,\ldots,R_k)\in \bbL[\tau(T')]$, 
    where $R_1,\ldots,R_k$ are the non-logical symbols occurring in $\psi$, we have
$$
\psi(R_1,\ldots,R_k)\in \tT'\quad\Leftrightarrow\quad \psi(\bfh(R_1),
\ldots,\bfh(R_k))\in \tT''.
$$

\end{enumerate}

Now  note that 

$\boxplus_1$  
 $\bbT/{\cong}$ has cardinality $\leq
2^{|\tau_0|+|\eps|+\theta}$.
\medskip 

Now let $\{\tT'_\alpha:\alpha<2^{|\tau_0|+|\eps|+\theta}\}$ be 
a list of members of $\bbT$ such that every isomorphism equivalence class 
over ${\rm Th}(M_1)$ is represented, and 
$\LL \tau(\tT'_\alpha)\setminus\tau_1:\alpha<2^{|\tau_0|+
|\eps|+\theta}\RR$ are pairwise disjoint.

Note  that 
${\rm Th}(M_1)\subseteq \tT'_\alpha$. Let $\tT_2'=\bigcup
\{\tT_\alpha':\alpha< 2^{|\tau_0|+|\eps|+\theta}\}$ and note 
\begin{enumerate}
    \item[$\boxplus_2$] $\tT'_2$ is consistent.
\end{enumerate}

\par \noindent
[Why? By Robinson consistency theorem.]
\medskip

Let $\tT_2$ be any completion of $\tT'_2$. So condition (a) holds;
proving (b) should be easy.

Let us prove {\bf (a)}$^+$ of the claim; this is really the proof that a theory $T$,
$|\tT|<\theta$, has a model in $2^{<\theta}$ universal for models
of $\tT$ of cardinality $\leq\theta$. We shall define by induction on
$\alpha<\theta$, a theory $\tT^2_\alpha$ such that:
\begin{enumerate}
\item[(A)] $\tT^2_0={\rm Th}(M_1)$,
\item[(B)] $\tT^2_\alpha$ a f.o.~theory,
\item[(C)] $\tau^2_\alpha=\tau(\tT^2_\alpha)$ has cardinality $\leq
2^{|\tau_0|+|\eps|+|\alpha|+\aleph_0}$,
\item[(D)] $\tT^2_\alpha$, $\tau^2_\alpha$ are increasing continuous in
$\alpha$,
\item[(E)] \underline{if}\ $\tau_1\subseteq\tau'\subseteq\tau^2_\alpha$,
$|\tau'|\leq |\tau^1|+|\alpha|$, $\tau'\subseteq\tau''$, $\tau''\cap
\tau^2_\alpha=\tau'$, $\tT^2_\alpha\rest\bbL_{\tau'}\subseteq \tT''
\subseteq\bbL[\tau'']$, $\tT''$  complete and $|\tau''\setminus\tau'|=
\{R\}$,\\
\underline{then} we can find $R'\in \tau^2_{\alpha+1}\setminus\tau^2_\alpha$
such that of the same arity.
$$
\tT''[\mbox{replacing }R\mbox{ by }R']\subseteq \tT^2_{\alpha+1}
$$
\end{enumerate}
There is no problem to carry out the induction, and $\bigcup\limits_{\alpha<
\theta}\tT^2_\alpha$ is as required.
\medskip

\noindent 2)\qquad Similar. 
\end{PROOF}

\begin{Explanation}
\label{1.7D}
{\rm

So for $M\models \tT$, $\bar b\in {}^{\kappa>}M$, we can choose 
$T_2[\bar b] \supseteq {\rm Th}(M,\bar b)$ depending on ${\rm Th}(M,\bar b)$ 
only, such that:
\begin{enumerate}
    \item[($\otimes$)] $M\models$ ``$\tT$ is $\kappa$--resplendent
    \underline{if}\ for every $\bar b\in {}^{\kappa>}M$, $(M,\bar b)$ is
    expandable to a model of $\tT_2[\bar b]$.
\end{enumerate}
W.l.o.g.~$\tau(\tT_2[\bar b])$ depends on $\lh(\bar b)$ and $\tau_0$
only, so it is $\tau(\lh(\bar{b}),\tau_M)$.

The things look quite finitary \underline{but} $\tT_2[\bar b]$ is not
continuous in ${\rm Th}(M,\bar b)$. I.e., $(\ast)\not\Rightarrow(\ast\ast)$,
where
\begin{enumerate}
\item[$(\ast)$] $\bar b^\alpha\in {}^{\kappa>}M$, for $\alpha\leq\delta$,
($\delta$ a limit ordinal) $\lh(\bar b^\alpha)=\eps$, and for every $n$,
$i_1<\ldots<i_n < \eps $  
 and a formula $\varphi(x_{i_1},\ldots,x_{i_n})\in
\bbL ( \tau _M )  $ 
for some $\beta<\delta$:
$$
\beta\leq\alpha\leq\delta\quad\Rightarrow\quad M\models\varphi[
b^\beta_{i_1},\ldots,b^\beta_{i_n}]\equiv\varphi[b^\alpha_{i_1},\ldots,
b^\alpha_{i_n}],
$$
\item[$(\ast\ast)$] for any $\varphi\in \bbL(\tau_2)$
for some $\beta<\delta$:
$$
\beta\leq\alpha\leq\delta\quad\Rightarrow\quad\big[\varphi\in \tT_2[
\bar{b}^\alpha]\ \Leftrightarrow\ \varphi\in \tT_2[\bar{b}^\beta]\big].
$$
\end{enumerate}
[You can make $\tT_1[\bar{b}]$ a somewhat continuous function of the
sequence $\bar b$ if we look at sub-sequences as approximations
rather than the type, but this is not used.]

For example, consider the case 
\begin{itemize}
    \item $M \models (\exists! {}^n x) R(x, b_n)$
    
    \item $\psi(y)$ says $F$ is a one-to-one function from $\{ x : R(x,y)\}$
    but not onto it.
\end{itemize}
}

This explains why you need ``infinitary Skolem functions''. 
\end{Explanation} 

\medskip

We shall try to construct $M$ such that for every $\bar b\in {}^\eps\! M$,
$(M,\bar b)$ is expandable to a model of $\tT_2[\bar b]$, so if
$\tau_2^\eps=\tau(\tT_2[\bar{b}])\setminus\tau(M,\bar b)$, this
means we have to define finitary relations/functions $R_{\bar b}$ (for $R\in
\tau_2^\eps$). We write here $\bar b$ as a sequence of parameters but
from another prospective the predicate/function symbol $R_{\square}(-)$ has
$\eps+{\rm arity}(R)$--places.

\begin{Explaining the first construction}
\label{1.7F}
{\rm
(i.e., \ref{1.28} below)\\
Eventually we build a generalization of ${\rm GEM}({}^{\kappa\geq}\lambda,\Psi)$, a
model with skeleton $\bar a_\eta$ ($\eta\in {}^{\kappa\geq}\lambda$)
witnessing $\kappa<\kappa(\tT)$, but the functions have any $\alpha<
\kappa$ places but not $\kappa$, and the indiscernibility 
demand is weak. We start as in \cite[\S2]{Sh:E59}, so for some formulas 
$\LL \varphi_\alpha(x, \bar y, \alpha) : \alpha < \kappa \RR$ we have
(where $\bar a_\eta=\LL a_\eta\RR$
for $\eta\in {}^{\kappa}\lambda$):
\[\eta\in {}^\kappa\lambda\ \&\ \nu\in {}^{\alpha+1}\lambda\quad\Rightarrow
\quad\varphi_\alpha(a_\eta,\bar a_\nu)^{{\rm if}(\nu\lhd
\eta)}.\]
Recall that for a formula $\varphi$ we let $\varphi^0 = \varphi$ and 
$\varphi^1 = \neg\varphi$, so $\varphi^{\mathrm{if}(\psi)}$ is either $\varphi$ 
or $\neg\varphi$ according to the truth value of $\psi$. Without loss of 
generality, for any $\alpha<\kappa$ for some sequence 
$\bar{G}_\alpha = \big\LL G_{\alpha,\ell} : \ell < \lh (\bar{y}_\alpha) \big\RR$
of unary function symbols such that for any $\eta \in {}^\kappa\lambda$, 
$$\bar a_{\eta\rest\alpha} = \bar{G}_\alpha (a_\eta) := \big\LL G_{\alpha,\ell} 
(a_\eta): {\ell} < \lh (\bar y_\alpha)\big\RR$$ 
so we can look at $\{a_\eta:\eta\in {}^\kappa\lambda\}$ as
generators. For $W\in [{}^\kappa\lambda]^{<\kappa}$, let $N_W=N[W]$ be the
submodel which $\{a_\eta:\eta\in W\}$ generates. So we would like to have:

\begin{enumerate}
\item[$(\alpha)$] $N_W$ has the finitary Skolem function (for $\tT$), and
moreover

$N_W$ has the finitary Skolem function for $\tT_2[\bar{b}]$ for each
$\bar{b}\in {}^{\kappa>}(N_W)$,
\item[$(\beta)$] monotonicity:\quad $W_1\subseteq W_2\ \Rightarrow\ N_{W_1}
\subseteq N_{W_2}$.
\end{enumerate}
So if $\cU\subseteq {}^\kappa\lambda$, then $N[\cU]=\{N_W:W\in
{}^{\kappa>}[\cU]\}$ is a $\kappa$--resplendent model of cardinality
$\lambda$.
\begin{enumerate}
\item[$(\gamma)$] Indiscernibility:\quad (We use here very ``minimal''
requirement (see below) but still enough for the omitting type in (1)
below):
\begin{enumerate}
\item[$(1)$] $\eta\in {}^\kappa\lambda\setminus \cU\quad\Rightarrow\quad
N[\cU]$ omits $p_\eta =:\{\varphi_\alpha(x,\bar a_{\eta\rest\alpha})
:\alpha<\kappa\}$;
(satisfaction defined in $N[^{\kappa }  
\lambda]$), 

\item[$(2)$] $\eta\in {}^\kappa\lambda\cap\cU\quad\Rightarrow\quad
N[\cU] $ 
realizes $ p_\eta$.
\end{enumerate}
\end{enumerate}
Now (2) was already guaranteed: $a_\eta$ realizes $p_\eta$.

For (1) it is enough
\begin{enumerate}
\item[$(1)'$] \underline{if}\ $W\in {}^{\kappa>}[{}^\kappa\lambda]$, $\eta\in
{}^\kappa\lambda\setminus W$ then $p_\eta$ is omitted by $N_W$ (satisfaction
defined in $N[^{\kappa 
}\lambda]$).
\end{enumerate}
Fix $W$, $\eta$ for $(1)'$. A sufficient condition is
\begin{enumerate}
\item[$(1)''$] for $\alpha<\kappa$ large enough, $\LL\bar a_{\eta
\rest\alpha\conc\LL i\RR}:i<\lambda\RR$ is indiscernible
over $N_W$
in $N[{}^\kappa\lambda]   $. 
\end{enumerate}
[if $\kappa(T)<\infty$, this immediately suffices; in the general 
case, and avoiding classification theory, use
\[
p'_\eta=\{\varphi_\alpha(x,\bar a_{\eta\rest(\alpha+1)})\ \&\ \lnot
\varphi(x,\bar a_{\eta\rest\alpha \caret \LL\eta(\alpha)+1\RR}):
\alpha<\kappa\}
\]
so
we use
\[\varphi' (x,\bar y'_\alpha)=\varphi_\alpha
(x,\bar y'_\alpha \rest 
 \lh(\bar y_\alpha)) \wedge 
\neg
\varphi_\alpha (x,\bar y'_\alpha\rest (\lh (\bar y_\alpha),
 2 
\lh (\bar y_\alpha)))\]
in the end].

Note: as $|W|<\kappa$, 
for some $\alpha(\ast)<\kappa$,
for every $\eta\in {}^{\kappa}\lambda$ 
\[W\cap\{\nu:\eta\rest\alpha(\ast)\lhd\nu\in {}^{\kappa}
\lambda\} \mbox{ is a singleton}\]
and $W\in\bfW_{\alpha(\ast)}$ (see below), this will be enough to omit the
type. The actual indiscernibility is somewhat stronger.
}
\end{Explaining the first construction}

\par \noindent
\underline{Further Explanation:} 
On the one hand, we would like to deal with
arbitrary
sequences of length $<\kappa$, on the other hand, we would like to retain
enough freedom to have the weak indiscernibility. What do we do? We define
our ``$\Phi$'' (not as nice as in \cite[\S2]{Sh:E59}, 
i.e., \cite[Ch.VII \S 3]{Sh:a}) by 
$\kappa$ approximations indexed for $\alpha\leq\kappa$.

For $\alpha\leq\kappa$, we essentially have $N_W$ for
\[\begin{array}{ll}
W\in\bfW_\alpha =: \{W:&W\subseteq {}^{\kappa 
}\lambda,\ |W|<\kappa \mbox{
and}\\
&\mbox{the function } 
\eta\mapsto\eta\rest\alpha\ (\eta\in W)\mbox{ is one-to-one}\}.
  \end{array}\]
Now, $\bfW_\alpha$ is partially ordered by $\subseteq$ but (for $\alpha < \kappa$) 
is not directed. For $\alpha<\beta$ we have $\bfW_\alpha\subseteq \bfW_\beta$ and $\bfW_\kappa=\bigcup\limits_{\alpha<\kappa}\bfW_\alpha$
is equal to $[{}^{\kappa}\lambda]^{<\kappa}$.

So if we succeed to carry out the induction for $\alpha<\kappa$, arriving at
$\alpha=\kappa$ the direct limit works and no new sequence of length
$<\kappa$ arises.

\newpage

\section{Proof of the Main Lemma}
In this section we get many models using a weak version of indiscernibility.

\begin{Hypothesis}
\label{1.17}
{\rm
\begin{enumerate}[(1)]
    \item $\tT$ is a fixed complete first order theory, $\kappa<\kappa(\tT)$, 
    $\bar{\varphi} = \LL\varphi_\alpha(x,y) : \alpha < \kappa \RR$ is a
    fixed witness for $\kappa<\kappa(\tT)$, recalling \ref{z17} and \ref{1.19}. This means:
    \begin{enumerate}
        \item[$(*)$] for any $\lambda$, for some model $M$ of $\tT$ and sequence
        $\LL a_\eta:\eta\in {}^{\kappa\geq}\lambda\RR$ with $a_\eta \in M$ we
        have: if $\eps<\kappa$, $\eta\in {}^\kappa\lambda$, $\alpha<\lambda$ then
        $$M \models \varphi_\eps[a_\eta, 
        a_{\eta\rest\eps\conc\LL\alpha\RR}]^{{\rm if}(\alpha=\eta(\eps))}.$$
    \end{enumerate}
    \item Let $\mu$ be infinite large enough cardinal; $\mu=\beth_\omega(|T|)$ is O.K.
\end{enumerate}
}
\end{Hypothesis}

\noindent{\sc Remark:}\label{1.18}
Why are we allowed in \ref{1.17}(1) 
to use $\varphi_\alpha(x,y)$ instead $\varphi(x,
\bar{y})$? We can work in $\tT^{{\rm eq}}$, see \cite[Ch.~III]{Sh:c}
 and anyhow this is, in fact,
just a notational change.
\medskip

\begin{Definition}
\label{1.19}
\begin{enumerate}
    \item For $\alpha<\kappa$ and $\rho\in {}^\alpha\mu $, let 
    $I_\rho = I^\alpha_\rho = I^{\alpha,\mu}_\rho$ be the model
\[
\big( \{\nu\in (^\kappa\mu):\nu\rest\alpha=\rho\},
E_i,
{<_i}\big)_{i<\kappa},
\]
    where
    \[
    \begin{array}{lcl}
    E_i&=&\left\lbrace(\eta,\nu):\eta\in
    {}^\kappa\mu,\nu\in {}^\kappa 
    \mu, \eta
    \rest i=\nu\rest i\right\rbrace,\\
    <_i&=&\left\lbrace(\eta,\nu):\eta\ E_i\nu\ \mbox{ and }\ \eta(i)<\nu(i)
    \right\rbrace.
    \end{array}
    \]
    
    \item Let $\bfW_\alpha=\bfW_\alpha^\mu=\{W\subseteq{}^\kappa\mu : W$ 
    has cardinality $<\kappa$ and for any $\eta\neq\nu$ from $W$ we have 
    $\eta\rest\alpha\neq\nu\rest\alpha\ \}$, and 
    $\bfW_{<\alpha} = \bigcup\limits_{\beta<\alpha}\bfW_\beta$.
    
    \item We say that $\bfW$ is $\alpha$--invariant, or 
    $(\alpha,\mu)$--invariant, when $\bfW\subseteq\bfW_\alpha$ and:

    \underline{if}\ $W_1,W_2\in\bfW_\alpha$, $h$ is a one-to-one function 
    from $W_1$ onto $W_2$ and $\eta\rest\alpha = h(\eta)\rest\alpha$ for
    $\eta\in\bfW_1$,

    \underline{then} $W_1\in\bfW\ \Leftrightarrow\ W_2\in\bfW$.
    
    \item We say $\bfW \subseteq \bfW_\alpha$ is hereditary if 
    $W'\subseteq W \in \bfW \Rightarrow W' \in \bfW$.
    
    \item let $\bfW_{<\alpha} = \bigcup\{ \bfW_\beta : \beta < \alpha\}$,
    so obviously it is contained in $\bfW_\alpha$. 
\end{enumerate}
\end{Definition}

\begin{Definition}
\label{1.20}
\begin{enumerate}
    \item Let $\theta = \theta_{T, \kappa}$ be the minimal cardinal satisfying:
    \begin{enumerate}
        \item[(a)] $\theta=\theta^{<\kappa}\geq |\tT|$,
        
        \item[(b)] if $M$ is a model of $\tT$, $\bar{b}\in {}^{\kappa>}M$,
        \underline{then} there is a complete (first order) theory $\tT^*$
        of cardinality $\leq \theta$ with Skolem functions extending 
        ${\rm Th}(M,\bar{b})$ such that:

        if $\tT'\supseteq {\rm Th}(M,\bar{b})$ and 
        $\tau(\tT')\setminus\tau_{(M,\bar{b})}$ has cardinality $<\kappa$
         \underline{then} there is a one-to-one mapping from $\tau(\tT')$ into 
        $\tau(\tT^*)$ over $\tau_{(M,\bar{b})}$ preserving arity and being a 
        predicate / function symbol, and mapping $\tT'$ into $\tT^*$.
    \end{enumerate}
    
\item For $\eps<\kappa$, let $\tau[\tT,\eps]$ be a
vocabulary consisting of $\tau_{T}$, the individual constants
$\underline{b}_\xi$ for $\xi<\eps$, and the $n$--place predicates
$R_{T,j,n}$ for $j<\theta$ and $n$--place function symbols $F_{T,j,n}$ for
$j<\theta$.

For $\eps<\kappa$ and a complete theory $T^\oplus$ in the
vocabulary $\tau_{T}\cup\{\underline{b}_\xi:\xi<\eps\}$ extending
$T$, let $T^*[T^\oplus]$ be a complete first order theory in the
vocabulary $\tau[T,\eps]$ such that if $(M,\bar{b})$ is a model
of $T^\oplus$, then $T^*[T^\oplus]$ is as in clause (b) of
part (1).
\item For $M\models T$ and $\eps<\kappa$ and $\bar{b}\in
{}^\eps M$, let $T^*[\bar{b},M]=T^*[{\rm Th}(M,\bar{b}
)]$.
\end{enumerate}
\end{Definition}

\begin{Remark}
\label{2.3A}
{\rm
Note that $\theta$ is well defined by \ref{1.7C}. 
 In fact, $\theta =
\Pi \{ 2^{ |T| + \sigma } : \sigma^+ < \kappa \} $  
is OK.
}
\end{Remark}

\begin{Main Definition}
\label{1.21}
We say that $\gm$ is an approximation (or an $\alpha$--approximation,
or $(\alpha,\mu)$--approximation) \underline{if}
\begin{enumerate}
\item[$(*)_1$] $\alpha\leq\kappa$ (so $\alpha=\alpha_\gm=\alpha(\gm)$),
\item[$(*)_2$] $\gm$ consists of the following (so we may give them
subscript or superscript $\gm$):
\begin{enumerate}
\item[(a)] a model $M=M_\gm$;
\item[(b)] a set $\cF=\cF_\gm$ of symbols of functions,
each $f\in \cF$ has an
interpretation, a function $f_\gm$ with range $\subseteq M$,
but when no confusion arises we may write $f$ instead of
$f_\gm$, (or $f^\gm$, note that the role of those
$f$-s is close to that of function symbols in vocabularies,
but not equal to);  
\item[(c)] each $f\in \cF$ has $\zeta_f<\kappa$ places, to each place
$\zeta$ (i.e., an ordinal $\zeta<\zeta_f$) a unique $\eta_\zeta\in {}^\alpha
\mu$, $\eta_\zeta=\eta^f_\zeta=\eta(f,\zeta)$ is attached such that
\[[\zeta\neq\xi\quad \Rightarrow\quad \eta_\zeta\neq \eta_\xi],\]
and the $\zeta$-th variable of $f$ varies on $I_{\eta_\zeta}$, i.e.,
$f_\gm(\ldots,x_\zeta,\ldots)_{\zeta<\zeta_f}$ is well defined iff
$\bigwedge\limits_{\zeta<\zeta_f}x_\zeta\in I_{\eta_\zeta}
=I^{\alpha,\mu}_{\eta_\zeta}$;

we may write $f_\gm(\ldots,\nu_\eta,\ldots)_{\eta\in w[f]}$ instead
$f_\gm(\ldots,\nu_{\eta(f,\zeta)},\ldots)_{\zeta<\zeta_f}$, where
$w[f] = \{\eta(f, \zeta) : \zeta < \zeta_f\}$; and 
$$f \in \cF \Rightarrow (\exists W \in \bfW) \big[w[f] =\{\eta\rest\alpha : \eta\in W\} \big]$$
see clause (e) below;
\item[(d)] for each $\bar b\in {}^{\kappa>}|M|$, an expansion $M_{\bar b}$ of
$(M, \bar{b})$ to a model of $T^*[\bar{b},M]$, (see above in Definition
\ref{1.20}; so $M_{\bar{b}}$ has Skolem functions and it witnesses
$\kappa$--resplendence for this sequence in $M$);
\item[(e)] $\bfW=\bfW_\gm\subseteq\bfW_\alpha$
which is $\alpha$--invariant and hereditary;
\item[(f)] for $W\in\bfW$, $N_W$ which is the submodel of
$M$ with universe
\[\begin{array}{lr}
\{f(\ldots,\eta_\zeta,\ldots)_{\zeta<\zeta_f}:&f\in \cF,\ f(\ldots,
\eta_\zeta,\ldots)_{\zeta<\zeta_f}\mbox{ well defined, }\ \\
&\mbox{and }\eta_\zeta\in W\mbox{ for every }\zeta\},
\end{array}\]
\item[(g)] a function $\bff= \bff_\gm$,
\end{enumerate}
\end{enumerate}
\underline{such that} $\gm$ satisfies the following:
\begin{enumerate}
    \item[(A)] $M$ is a model of $T$,
    
    \item[(B)] [\underline{witness for $\kappa<\kappa(T)$:}]\quad for our
    fixed sequence of first order formulas 
    $\LL \varphi_\zeta(x,y) : \zeta < \kappa \RR$ from $\bbL(\tau_T)$ 
    $($depending neither on $\alpha$ nor on $\gm)$ we have
    $f^\ast_{\rho,\zeta}\in \cF$ for $\zeta\leq\kappa$, $\rho\in {}^\alpha\mu$
    $($we also call them $f_{\rho,\zeta}^\gm)$ such that:
    \begin{enumerate}
        \item[(i)]   $f^\ast_{\rho,\zeta}$ is a one place function, with
        $\zeta_{f^*_{\rho,\zeta}}, \eta_0^{f^*_{\rho,\zeta}}$ from clause (c) being
        $1,\rho$ respectively.
        
        \item[(ii)]  $f^\ast_{\rho_1,\zeta}(\nu_1)=f^\ast_{\rho_2,\zeta}(\nu_2)$ if
        $\nu_1\rest\zeta=\nu_2\rest\zeta$ and they are well defined, i.e.
        $\rho_\ell \triangleleft \nu_\ell \in {}^\kappa\mu$,
        
        \item[(iii)] if $\rho_\ell\in {}^\alpha\mu$, $\nu_\ell\in I^\alpha_{\rho_\ell}$ 
        for $\ell=1,2$ and $\zeta<\kappa$ \underline{then}:
        \[
        M \models\varphi_\zeta\left[ f^\ast_{\rho_1,\kappa}(\nu_1),
        f^\ast_{\rho_2,\zeta+1}(\nu_2)\right]\qquad\mbox{ iff }\qquad[\nu_1
        \rest(\zeta+1)=\nu_2\rest(\zeta+1)],
        \]
    \end{enumerate}
    
    \item[(C)] $N_W\prec M$ for $W\in \bfW$; moreover, $N_W = M_\gm \rest A_W$, where $A_W$ 
    is the minimal subset of $M_\gm$ such that: (see below)
    
    \noindent if $\bar f \in \cF$, $\bar \nu \subseteq W$, $\bar \nu \in \dom(\bar f)$, 
    \underline{then} $|N_{\bar b}| \subseteq A_W$. Therefore (see \ref{1.22}) 
    if $\gm$ is full and closed under terms then $M_{\bar b} \rest A_W \prec M_{\bar b}$.

    \item[(D)] [\underline{$\bff=\bff_\gm$ witnesses an amount of resplendence}]
    \begin{enumerate}
        \item[$(\alpha)$] the domain of $\bff$ is a \underline{subset} of $\bbF_\gm$, 
        where $\bbF_\gm$ is the subset of 
        $$\{\bar f = \LL f_\eps : \eps < \eps_{\bar f}\RR : \eps_{\bar f} < \kappa,\ f_\eps \in \cF\}$$ 
        for which $\zeta_{f_\eps} =: \zeta_{\bar f}$ does not depend on $\eps$, the sequence 
        $\eta(f_\eps, \zeta) =: \eta(\bar f, \zeta)$ does not depend on $\eps$ for all 
        $\zeta < \zeta_{\bar f}$, and $\{\eta(\bar f, \zeta) : \zeta < \zeta_{\bar f}\} \subseteq W$
        for some $W \in \bfW_\gm$.
        
        [so $\bar f \in \bbF_\gm$ maps $\prod\limits_{\zeta < \zeta_{f}} I_{\eta(\bar f, \zeta)}$ into ${}^{\eps(f)}\!M$]
        
        \item[$(\beta)$] for $\bar{f}\in\Dom(\bff)$, $\bff(\bar{f})$ is a
        function with domain
        \[
        \begin{array}{ll}
        \{\sigma(\bar{x}):&\sigma(\bar{x})
        \mbox{ is a $\tau[T, \eps_{\bar{f}}]$--term},\mbox{ and }
        \bar{x}=\LL x_\xi:\xi\in u\RR\\
        &\mbox{for some  finite subset }u=u_\sigma\mbox{ of }\eps_{\bar{f}}\;\}
        \end{array}
        \]
        and if $\sigma(\bar{x})\in\Dom(\bff(\bar{f}))$ \underline{then}
        \[
        \hspace{-1cm}\bff(\bar{f})(\sigma(\bar{x}))\in \cF[\bar f] := 
        \Big\{f\in\cF:\zeta_f = \zeta_{\bar{f}},\ 
        (\forall\zeta<\zeta_f) \big[ \eta(f,\zeta)=\eta(\bar{f},\zeta) \big] \Big\},
        \]
        
        \item[$(\gamma)$] 
        \begin{enumerate}
            \item[$\bullet$]  if $ \bar f \in \Dom(\bff)$ and 
            $\bar b = \big\LL f_\eps(\ldots, \nu_{\eta(\bar f, \zeta)},\ldots)_{\zeta < \zeta_{\bar f}} :
            \eps < \eps_{\bar f} \big\RR$ 
            \underline{then} the set
            $\big\{ f(\ldots, \nu_{\eta( \bar f, \zeta)}, \ldots)_{\zeta < \zeta _{\bar f}} : 
            f \in \cF[\bar f] \big\}$ is the universe of an elementary submodel of $M_{\bar b}$ 
            called $N_{\gm,\bar f}$. 
            {This actually follows from the next point:}
            
            \item[$\bullet$]  
            if $(\bff(\bar{f}))(\sigma(\bar{x}))=f^* \in \cF[\bar f]$,
            $W\in\bfW$, $\nu_\zeta\in W$ and $\nu_\zeta\rest\alpha = \eta(\bar{f},\zeta)$ 
            for $\zeta < \zeta_{\bar f}$, and $\bar x = \LL x_\xi : \xi \in u\RR$ (where 
            $u \subseteq \zeta_{\bar f}$),
            and $\bar{b} = \bar{f}(\bar{\nu}) = \LL f_\eps(\bar{\nu}) : \eps < \eps_{\bar{f}}\RR$, 
            \underline{then} 
            \[
            \sigma^{M_{\bar{b}}}(\LL f_\xi(\ldots,\nu_\zeta,\ldots)_{\zeta<\zeta_{
            \bar{f}}}:\xi\in u\RR)=f^*(\ldots,\nu_\zeta,\ldots)_{\zeta<\zeta_{\bar{f}}}.
            \]
        \end{enumerate}
    \end{enumerate}
\end{enumerate}

[explaining $(\gamma)$: we may consider $\bar{b} \in N_{W_1} \cap N_{W_2}$,
and we better have that the witnesses for resplendence demands,
specialized to
$\bar{b}$, in $N_{W_1}$ and in $N_{W_2}$
are compatible so that in the end
resplendence holds].
\end{Main Definition}

However, we shall not get far without at least more closure and coherence of
the parts of $\gm$.

\begin{Definition}
\label{1.22}
\begin{enumerate}[(1)]
\item An approximation $\gm$ is called full \underline{if}\
$\bfW_\gm=
\bfW_{\alpha(\gm)}$, and is called semi--full if
$\bfW_{<\alpha(\gm)}\subseteq\bfW_{\gm}\subseteq\bfW_{\alpha(\gm)}$ and is called
almost full if it is semi full when $\alpha$ is limit ordinal and full when
$\alpha$
is a non-limit ordinal.
\item An approximation $\gm$ is $\beta$--resplendent if $\beta\leq
\alpha_\gm$ and: (recalling $\dom(\bff_\gm) \subseteq \bbF_\gm$)
\[\begin{array}{ll}
\mbox{if }&W\in\bfW_\beta\cap\bfW_\gm
\mbox{ and }\bar{f}\in\bbF_\gm,\mbox{ and}\\
&\{\eta(\bar{f},\zeta):\zeta<\zeta_{\bar{f}}\}\subseteq\{\nu\rest
\alpha:\nu\in W\},\\
\mbox{then}&\bar{f}\in\Dom(\bff_\gm).
\end{array}\]
\item In part (2), if we omit $\beta$, we mean $\beta=\alpha_\gm$, and
``$<\beta^*$'' means for every $\beta<\beta^*$.
\item An approximation $\gm$ is called term closed \underline{if}:
\begin{enumerate}
\item[(E)] \underline{Closure under terms of $\tau$:}

Assume that $u\subseteq{}^\alpha\mu$, $|u|<\kappa$, and for some $W\in
\bfW_\gm$, $u\subseteq\{\eta\rest\alpha:\alpha\in W\}$,
and
$\LL\eta_\zeta:\zeta<\zeta^*\RR$ lists $u$ with no
repetitions and
$f_\ell\in \cF_\gm$, $\ell<n$, satisfies $\left\lbrace\eta
(f_\ell,\zeta):\zeta<\zeta_{f_\ell}\right\rbrace\subseteq u$,
$\sigma$ is an
$n$-place $\tau(T)$-term so $\sigma=\sigma (x_0,\ldots,x_{n-1})$.
\underline{then} for some $f\in \cF_\gm$ satisfying $\zeta_f=\zeta^*$,
$\eta(f,\zeta)=\eta_\zeta$, for any choice of
$\LL \nu_\eta:\eta\in u\RR$ such that $\eta\triangleleft
\nu_\eta\in {}^\kappa\mu$ for $\eta\in u$,
 and 
$\{\nu_\eta:\eta\in  u 
\}\subseteq W'\in\bfW$
 for some $W'$ we have
$$
f_\gm\left(\ldots,\nu_\eta,\ldots\right)_{\eta\in w[f]} = 
\sigma\left(\ldots, f^\gm_\ell(\ldots,\nu_\eta,\ldots)_{\eta\in w[f_\ell]},
\ldots\right)_{\ell<n}
$$
(This clause may be empty, but it helps to understand clause {\bf (F)};
note that it is not covered by \ref{1.21}(D)($\beta$) as the functions 
do not necessarily have the same domain, hence this says something even for
$\sigma$ the identity. This implies that in clause $(f)$ of Definition 
\ref{1.21} we can demand $\{\eta_\gm(f,\zeta) : \zeta < \zeta_f\} = W$.
In other words, in \ref{1.21}(D), given one $\bar f \in \dom(\bff_\gm)$
we can find others; here we claim the existence of $\bar f$ for a given
$\LL \eta(f, \zeta) : \zeta < \zeta_{\bar f} \RR$.)

\item[(F)] \underline{Closure under terms of $\tau(M_{\bar b})$:}

Assume that $u\subseteq{}^\alpha\mu$, $|u|<\kappa$,
and $\LL\eta_\zeta:
\zeta<\zeta^*\RR$ lists $u$ with no repetitions, and for some
$W\in\bfW_\gm$, $u\subseteq\{\eta\rest\alpha:\eta\in
W\}$. \underline{If}\ $n<\omega$ and $f^\ell\in \cF_\gm$ for $\ell<n$,
$\bar{f}=\LL f_\eps:\eps<\eps(*)\RR\in
\Dom (\bff_\gm)$, and
\[\begin{array}{rl}
\eta(f_\eps,\zeta)\in u&\mbox{ for }\zeta<\zeta_{f_\eps},\mbox{
and}\\
\eta(f^\ell,\zeta)\in u&\mbox{ for }\zeta<\zeta_{f^\ell},\mbox{ and}\\
b_\eps=f_\eps(\ldots,\nu_{\eta(f_\eps,\zeta)},\ldots)_\zeta&
\mbox{ for }\eps<\eps(\ast),
\end{array}\]
$\bar b=\LL b_\eps:\eps<\eps(\ast)\RR$ and $\sigma(x_0,
\ldots,x_{n-1})$ is a $\tau(M_{\bar b})$--term,\\
\underline{then} for some $f\in \cF_\gm$ we have $w[f]=u$ and:

if $\nu_\eta\in I^\alpha_\eta$ for $\eta\in u$ and
$\{\nu_\eta:\eta\in u\}
\in \bfW_\gm$, then
$$
\hspace{-.7cm}
f\left(\ldots,\nu_\eta,\ldots\right)_{\eta\in u} =
\sigma^{M_{\bar b}}\left(\ldots,f_\gm^\ell(
\ldots,\nu_{\eta(f^\ell,\zeta)},\ldots)_{\zeta<\zeta(f^\ell)},
\ldots\right)_{\ell<n}.
$$
\end{enumerate}
\end{enumerate}
\end{Definition}

\begin{observation}
\label{2.6A}
In Definition  
\ref{1.22}(4) 
in clauses (E),(F) it
suffice to restrict ourselves to the case 
$ n=1 $ and $ \sigma $ 
is the identity.
\end{observation}

\begin{PROOF}{\ref{2.6A}}
By \ref{1.21}(D)$(\gamma $). 
\end{PROOF} 

Of course some form of indiscernibility will be needed.
\medskip

\begin{Definition}
\label{1.22x}
\begin{enumerate}
\item[(1)] Let $\bbE$ be the family of equivalence relations $\bfE$ on
\[
\{\bar{\nu}\in {}^{\kappa>}\!\left({}^\kappa\mu\right) : \bar{\nu}
\mbox{ without repetitions}\},
\]
or a subset of it, such that
\[\bar{\nu}^1\mathrel{\bfE}\bar{\nu}^2\quad\Rightarrow\quad\lh(\bar{\nu}^1)
=\lh(\bar{\nu}^2).\]
\item[(2)] Let $\bbE_\alpha$ be the family of $\bfE\in\bbE$ such that
\[\bar{\nu}\in\Dom(\bfE)\quad\Rightarrow\quad\LL\nu_\zeta\rest
\alpha:\zeta<\lh(\bar{\nu})\RR\mbox{ is without repetitions.}\]
\item[(3)] Let $\bfE^0_\alpha\in \bbE_\alpha$ be the following equivalence
relation:
\[\begin{array}{llll}
\bar{\nu}^1\;\bfE^0_\alpha\;\bar{\nu}^2&\mbox{iff}&&\mbox{for some }
\zeta<\kappa\mbox{ we have}\\
&&\mbox{(i)}&\bar{\nu}^1,\bar{\nu}^2\in{}^\zeta\!\left({}^\kappa\mu\right),\\
&&\mbox{(ii)}&\nu^1_\eps\rest\alpha=\nu^2_\eps
\rest\alpha \mbox{ for }\eps<\zeta,\\
&&\mbox{(iii)}&\LL\nu^1_\eps\rest\alpha:\eps<\zeta\RR\mbox{
is with no repetitions},\\
&&\mbox{(iv)}&\mbox{the set }\{\eps<\zeta:\nu^1_\eps\neq
  \nu^2_\eps\}\mbox{ is finite.}
  \end{array}\]

\item[(3A)] We say that $(\bar{\nu}^1,\bar{\nu}^2)$ are immediate neighbours \underline{if}\
$\lh
(\bar{\nu}^1)=\lh (\bar{\nu}^2)$, and for some
$\xi<\lh (\bar{\nu}^1)$ we have $(\forall\eps<\zeta)(\eps\neq\xi\
\Leftrightarrow\ \nu^1_\eps= \nu^2_\eps)$;
so the difference with (3)(iv) is that ``finite'' is replaced by
``a singleton''.
\item[(4)] Let $\bfE^0_{{<}\alpha}$ be defined like $\bfE^0_\alpha$
strengthening clause (iii) to
\begin{enumerate}
\item[(iii)$^+$] for some $\beta<\alpha$, the sequence $\LL \nu^\ell_\eps
\rest\beta:\eps<\zeta\RR$ is with no repetitions.
\end{enumerate}
\item[(5)] For $\alpha<\kappa$ and 
$\bfW\subseteq\bfW_\alpha$ 
let 
\[
\begin{array}{ll}
\seq_\alpha(\bfW) = \big\{\bar{\nu} : &\bar{\nu} \in 
{}^{\kappa{>}}\!\left({}^\kappa\!\mu\right)
\mbox{ is with no repetitions,}\\   
&\mbox{and for some $W\in\bfW$ we have}\\
&\{\nu_\xi : \xi < \lh(\bar{\nu})\} \subseteq W,\mbox{ and hence}\\
&\LL\nu_\zeta\rest\alpha:\zeta<\lh (\bar \nu)
\RR\mbox{ is with no 
repetitions}\big\}. 
\end{array}
\]
\item[(6)] We define $\bfE^1_\alpha$ as we define $\bfE^0_\alpha$ in part (3)
above, omitting clause (iv). We define $\bfE^1_{<\alpha}$ analogously to (4).
\end{enumerate}
\end{Definition}

\noindent{\sc Remark:}
The reader may concentrate on $\bfE^0_\alpha$, so the
``weakly'' version below.
\medskip

\begin{Definition}
\label{1.23}
\begin{enumerate}[(1)]
\item An approximation $\gm$ is called $\bfE$--indiscernible \underline{if}
\begin{enumerate}
\item[(a)] $\bfE\in \bbE_{\alpha(\gm)}$ refine
$\bfE^1_{\alpha(\gm)}$,
\item[(b)] if $\bar{\nu}^1,\bar{\nu}^2\in \seq_{\alpha(\gm)}
(\bfW_{\gm})$ and $\bar{\nu}^1\;\bfE\;\bar{\nu}^2$,
\underline{then} there
is $g$ (in fact, a unique $g=g^\gm_{\bar{\nu}^1,\bar{\nu}^2})$
such that
\begin{enumerate}
\item[$(\alpha)$] $g$ is an $(M_\gm,M_\gm)$--elementary mapping,
\item[$(\beta)$]  $\Dom(g)=\{f(\LL\nu^1_{h(\zeta)}:\zeta<\zeta_f
\RR): f\in \cF_\gm$ and $h$ is a one-to-one function
from $\zeta_f$ into $\lh(\bar{\nu}^\ell)$ such that 
$\eta(f,\zeta) \lhd \nu^1_\zeta \in {}^\kappa\!\mu\}$,
\item[$(\gamma)$] $g(f(\LL\nu^1_{h(\zeta)}:\zeta<\zeta_f\RR))=
f(\LL\nu^2_{h(\zeta)}:\zeta<\zeta_f\RR)$ for $f,h$ as above;
\end{enumerate}
\item[(c)] Assume $\bar{\nu}^1,\bar{\nu}^2\in\seq_{\alpha(\gm)}(
\bfW_\gm)$, $\bar{f}^1,\bar{f}^2\in\Dom(\bff_\gm)$,
$\zeta^*=\zeta_{\bar{f}^1} =\zeta_{\bar{f}^2}$, and for some
one-to-one function $h$
from $\zeta^*$ to $\lh(\bar{\nu}^\ell)$ we have
$\eta(\bar{f}^\ell,\zeta)=\nu^m_{h(\zeta)} \rest \alpha$ for
$\ell,m=1,2$,
and $\bar{\nu}^1\;\bfE\;\bar{\nu}^2$. Let
\[
\bar{b}^\ell=\big\LL f^\ell_\eps\big(\LL\nu^\ell_{h(\zeta)} : 
\zeta < \zeta^*\RR \big) : \eps < \lh(\bar{f}^\ell)\big\RR.
\]
\underline{then} there is $g$ such that
\begin{enumerate}
\item[$(\alpha)$] $g$ is an $(M^\gm_{\bar{b}_1},M^\gm_{
\bar{b}_2})$--elementary mapping,
\item[$(\beta)$] $g=g^\gm_{\bar{\nu}^1,\bar{\nu}^2}$ from clause (b)
above.
\end{enumerate}
\end{enumerate}
\item An approximation $\mathfrak m$ is strongly indiscernible if it is ${\bf
E}^1_{\alpha(\gm)}$--indiscernible.
\item
\begin{enumerate}
\item[(a)] An approximation $\gm$ is weakly indiscernible when it is
$\bfE^0_{\alpha(\gm)}$-indiscernibility.
\item[(b)] An approximation $\gm$ is weakly/strongly nice
\underline{if}\ it is term closed and weakly/strongly indiscernible.
\item[(c)] An approximation $\gm$ weakly/strongly good \underline{if}\
it is weakly/strongly nice and is almost full.
\item[(d)] An approximation $\gm$ is weakly/strongly excellent
\underline{if}\ it is weakly/strongly good, and is resplendent,
see Definition \ref{1.22}(2),(3). 
\end{enumerate}
\end{enumerate}
\end{Definition}

\begin{Discussion}
\label{1/23k} 
Why do we have the weak and strong version?

In the proof of the main subclaim \ref{1.28} below 
the proof for the
weak
version is easier \underline{but} we get from it a weaker conclusion:
$\geq \lambda^+$ non-isomorphic $\kappa$-resplendent of cardinality
$\lambda=\lambda^\kappa$, 
whereas from the strong version we 
would 
get $2^\lambda$. 
But see \S3.
\end{Discussion}

\begin{Claim}
\label{1.24}
Let $\gm$ be an approximation.
\begin{enumerate}[(1)]
\item In the definition of ``$\gm$ is
$\bfE^0_\alpha$--indiscernible'', it is enough to deal with immediate
$\bfE^0_\alpha$--neighbors (see Definition \ref{1.22x}(3)).
\item If $\gm$ is weakly/strongly excellent \underline{then} $\gm$ is
weakly/strongly good.
\item If $\gm$ is $weakly/strongly$ good \underline{then} 
$\gm$ is $weakly/strongly$ nice.
\item If $\alpha_\gm=0, \bfE\in \bbE_\gm$ \underline{then}
$\gm$ is $\bfE$-indiscernible \underline{iff}\ $\gm$ is strongly indiscernible.
\end{enumerate}
\end{Claim}

\begin{Definition}
\label{1.22y}
\begin{enumerate}[(1)]
    \item For approximations $\gm_1,\gm_2$ let ``$\gm_1\leq_\bfh\gm_2$'' 
    or ``$\gm_1\leq \gm_2$ as witnessed by $\bfh$'' mean that:
    \begin{enumerate}
        \item[(a)] $\alpha(\gm_1)\leq\alpha(\gm_2)$,
        
        \item[(b)] $\bfW_{\gm_1}\subseteq \bfW_{\gm_2}$,
        
        \item[(c)] $\bfh$ is a partial function from $\cF_{\gm_2}$ into $\cF_{\gm_1}$,
        
        \item[(d)] if $\bfh(f_2)=f_1$ \underline{then} they have the same arity 
        (i.e., $\zeta^{\gm_1}_{f_1}=\zeta^{\gm_2}_{f_2}$) and
        \[
        \zeta<\zeta^{\gm_1}_{f_1}\quad\Rightarrow\quad\eta_{\gm_1}(
        f_1,\zeta)=\eta_{\gm_2}(f_2,\zeta)\rest\alpha(\gm_1),
        \]
        
        \item[(e)] if $f_1\in \cF_{\gm_1}$ and $W\in\bfW_{\gm_1}$ and
        \[
        \{\nu\rest\alpha(\gm_1):\nu\in W\}=\{\eta_{\gm_1}(f_1,\zeta) : 
        \zeta<\zeta^{\gm_1}_{f_1}\},
        \]
        \underline{then} there is one and only one $f_2\in \cF_{\gm_2}$
        satisfying
        \[
        \bfh(f_2) = f_1\quad\mbox{and}\quad\{\eta_{\gm_1}(f_2,\zeta) : 
        \zeta < \zeta^{\gm_1}_{f_2}\} = \{\nu\rest\alpha(\gm_2) : \nu \in W\},
        \]
        
        \item[(f)] for $W\in\bfW_{\gm_1}$, the mapping $g^{\gm_2}_{\gm_1}[W,\bfh]$ 
        defined below is an elementary embedding from $N^{\gm_1}_W$ into 
        $N^{\gm_2}_W$, where:
        \begin{enumerate}
            \item[$(*)$] if $f_1\in\cF_{\gm_1}$, $f_2\in \cF_{\gm_2}$ are as in 
            clause (e) (so $\bfh(f_2) = f_1$), and
            \[
            a = f_1^{\gm_1}(\ldots,\nu_\zeta,\ldots)_{\zeta<\zeta_{f_1}},\quad
            \mbox{ and }\quad\{\nu_\zeta:\zeta<\zeta^{\gm_1}_f\}\subseteq W
            \]
            (so $a\in N^{\gm_1}_W$), then $(g^{\gm_2}_{\gm_1}[W,\bfh])(a) = 
            f^{\gm_2}_2 (\ldots,\nu_\zeta,\ldots)_{\zeta<\zeta_{f_2}}$,
        \end{enumerate}
        
        \item[(g)] if $\bar{f}^1 = \LL f^1_\xi : \xi < \eps\RR \in \Dom(\bff_{\gm_1})$ 
        and $\eta_{\gm_1}(\bar{f}^1,\zeta) \unlhd \eta_\zeta\in {}^{\alpha(\gm_2)}\mu$ 
        for $\zeta < \zeta^\gm_{\bar{f}^1}$, 
        $\bar{f}^2 = \LL f^2_\xi : \xi < \eps\RR \in {}^\eps(\cF_{\gm_2})$, and
        $\zeta^{\gm_2}_{\bar{f}^2}=\zeta^{\gm_1}_{\bar{f}^1}$, and
        $$\xi<\eps \wedge \zeta<\zeta_{\bar{f}^1}\ \Rightarrow\ 
        \eta(f^2_\xi,\zeta)=\eta_\zeta \wedge \bfh(f^2_\xi) = f^1_\xi$$ 
        \underline{then}
        \begin{enumerate}
            \item[$(\alpha)$] $\bar{f}^2\in\Dom(\bff_{\gm_2})$,
            
            \item[$(\beta)$]  $\bfh\left((\bff_{\gm_2}(\bar{f}^2))(\sigma(\LL x_\xi : 
            \xi \in u \RR)) \right) = \big( \bff_{\gm_1}(\bar{f}^1)\big)
            (\sigma(\LL x_\xi : \xi \in u\RR))$, when $u$ is a finite subset of $\eps$
            
            \item[$(\gamma)$] Assume $\nu_\zeta\in I_{\eta_\zeta}$ for 
            $\zeta < \zeta^{\gm_1}_{\bar{f}^1}$, and 
            $W = \{\nu_\zeta : \zeta < \zeta^{\gm_1}_{\bar{f}^1}\}$, 
            $\bar{b}^\ell = \LL f^\ell_\xi(\ldots, \nu_\zeta,\ldots)_\zeta : \xi < \eps\RR$.
            Then the mapping $g^{\gm_2}_{\gm_1}[W,h]$ (see clause (f) above) is a
            $(M_{\bar b_1},M_{\bar b_2})$-elementary mapping from 
            $M^{\gm_1}_{\bar{b}^1}\rest |N^{\gm_1}_W|$ 
            onto $M^{\gm_2}_{\bar{b}^2}\rest |N^{\gm_2}_W|$.
        \end{enumerate}
    \end{enumerate}
    
    \item We say that $\LL \gm_\beta, \bfh^\beta_\gamma : \beta < \alpha, \gamma\leq\beta\RR$ 
    is an inverse system of approximations \underline{if}
    \begin{enumerate}
        \item[(a)] $\gm_\beta$ is a $\beta$--approximation (for $\beta<\alpha$),
    
        \item[(b)] $\gm_\gamma\leq_{\bfh^\beta_\gamma} \gm_\beta$ for $\gamma\leq\beta$,
    
        \item[(c)] $\bfh^\beta_\beta$ is the identity,
        
        \item[(d)] if $\beta_0<\beta_1<\beta_2<\alpha$ then $\bfh^{\beta_2}_{\beta_0} = \bfh^{\beta_1}_{\beta_0} \circ \bfh^{\beta_2}_{\beta_1}$.
    \end{enumerate}
    
    \item We say that an inverse system of approximations 
    $\LL \gm_\beta, \bfh^\beta_\gamma : \beta < \alpha, \gamma \leq \beta \RR$ 
    is continuous at $\delta$ \underline{if}:
    \begin{enumerate}
        \item[(a)] $\delta<\alpha$ is a limit ordinal,
        
        \item[(b)] $\bfW_{\gm_\delta}=\bigcup \{\bfW_{\gm_\beta} : \beta < \delta\}$,
        
        \item[(c)] $\cF_{\gm_\delta} = \bigcup\{\Dom(\bfh^\delta_\beta) : \beta < \delta\}$,
        
        \item[(d)] $\Dom(\bff_{\gm_\delta})=\{\bar{f}^2:$ for some $\beta<\delta$ and 
        $\bar{f}^1 \in \Dom(\bff_{\gm_1})$ of length $\lh(\bar{f}^2)$ we have 
        $\bfh^\delta_\beta(f^2_\xi) = f^1_\xi$ for $\xi < \lh(\bar{f}^2)\}$.
    \end{enumerate}
\end{enumerate}
\end{Definition}

\begin{Discussion}\label{b35}
Having chosen above our order, when can we get the appropriate indiscernibility?
As we are using finitary partition theorem (with finitely many colours), 
we cannot make the type of candidates for fixed $\bar b$.
However, we may have, \emph{a priori}, enough indiscernibility to fix 
the type of enough $\bar b'$-s and then use the existence of indiscernibles  
to uniformize the related $M_{\bar b}$-s.
\end{Discussion}

\begin{Claim}\label{1.23x}
1) There is an excellent 0--approximation.

2) Moreover, there is a $\Phi$ such that:
\begin{enumerate}[(a)]
    \item $\Phi$ is a template proper for the tree $I_{\LL\ \RR}$,
    
    \item $\Phi$ is nice (see \cite[1.8]{Sh:E59}), 
    that for $ {\ell} $ the logic $ \mathbb{L} $, first order logic, the default value,
    
    \item $M_* = \mathrm{GEM}(I_{\LL\ \RR}, \Phi)$ is a model of $T$,
    
    \item if $\eta \in {}^\kappa\mu$ then $\mathrm{GEM}(\{\eta\}, \Phi)$ 
    is $\kappa$-resplendent,
    
    \item $\LL a_\eta : \eta \in {}^{\kappa\geq}\!\mu \RR$, 
    $\LL F_\eps^{M_*} : \eps < \kappa \RR$ are as above. 
\end{enumerate}
\end{Claim}

\begin{PROOF}{\ref{1.23x}} 
Recall that the sequence $\LL\varphi_\alpha(x,y) : \alpha < \kappa\RR$
exemplifies $\kappa<\kappa(T)$, see \ref{1.17} above. Hence by
clause (b) of 
     \cite[1.10(3)]{Sh:E59},   
we can find a template $\Phi$ proper for the tree $I_{\LL\ \RR}$, i.e.,
${}^{\kappa\geq}\mu$, with skeleton $\LL a_\eta:\eta\in{}^{\kappa\geq}
\mu\RR$ such that for $\nu\in{}^\kappa\mu$ and $\rho\in {}^{\alpha+1}
\mu$ we have
\[
{\rm GEM} 
({}^{\kappa\geq}\mu,\Phi) \models \varphi_\alpha(a_\nu,a_\rho)\quad
\mbox{ iff }\quad\rho \lhd \nu.\]
Without loss of generality, for some unary function symbols $F^*_\eps
\in\tau(\Phi)$, we have
${\rm GEM}({}^{\kappa\geq}\mu,\Phi)\models\mbox{``}F_\eps(
a_\eta)=a_{\eta\rest\eps}\mbox{''}$ for $\eta\in {}^\kappa
\mu$. Now, by induction on $\eps<\kappa$ we choose $\Phi_\eps$
such that:
\begin{enumerate}
    \item[(a)] $\Phi_\eps$ is a template proper for ${}^{\kappa\geq}\mu$
    which is nice (see 
     \cite[1.7]{Sh:E59} + 
     \cite[1.8(3),(4)
     ]{Sh:E59}), 
     
    \item[(b)] $\tau(\Phi_\eps)$ has cardinality $\leq\theta$ (see Definition \ref{1.20}),
    
    \item[(c)] $\Phi_0=\Phi$,
    
    \item[(d)] the sequence $\LL\Phi_\eps:\eps<\kappa\RR$
    is increasing with $\eps$, that is  (see \cite[1.8(1B)]{Sh:E59}: 
    \[
    \zeta<\eps\quad\Rightarrow\quad \tau(\Phi_\zeta)\subseteq\tau(
    \Phi_\eps)\ \mbox{ and }\ {\rm GEM}_{\tau (T)}({}^{\kappa\geq}\mu, \Phi_\zeta)\prec
    {\rm GEM}_{\tau (T)}({}^{\kappa\geq}\mu,\Phi_\eps),
    \]
    
    \item[(e)] the sequence $\LL\Phi_\eps : \eps < \kappa\RR$
    is continuous, i.e., if $\eps$ is a limit ordinal then $\tau(
    \Phi_\eps)=\bigcup\limits_{\zeta<\eps}\tau(\Phi_\zeta)$,
    
    \item[(f)] if $\bar{\sigma}=\LL\sigma_i(x):i<i^*\RR$ is a sequence
    of length $<\kappa$ of unary terms in $\tau(\Phi_\eps)$, and
    $M^{\eps+1}={\rm GEM}_{\tau (T)}({}^{\kappa\geq}\mu,\Phi_{\eps+1})$, and
    for $\nu\in {}^\kappa\mu$ we define $\bar{b}=\bar{b}_{\bar{\sigma},\nu}$ as
    \[
    \LL\sigma^{M^{\eps+1}}_i(a_\nu):i<i^*\RR\in {}^{i^*}
    ({\rm GEM}_{\tau (T)}(\{\nu\},\Phi_{\eps+1})),
    \]
    \underline{then} we can interpret a model $M^{\eps+1}_{\bar{b}}$ of
    $T^*[\bar{b},M^{\eps+1}\rest\tau_{T}]$ in $M^{\eps+1}$, which means:
    \begin{enumerate}
        \item[$(\alpha)$] if 
        $R\in \tau_{T}[\bar{b},M^{\eps+1}\rest\tau(T)]\setminus\tau_T$ 
        is a $k$--place predicate, then there is a $(k+1)$--place predicate $R_*\in\tau(\Phi_{\eps+1})\setminus\tau(\Phi_\eps)$ such that
        \[
        M^{\eps+1}_{\bar{b}}\models R[c_0,\ldots, c_{k-1}]\quad\mbox{ iff }
        \quad M^{\eps+1}\models R_*[c_0,\ldots,c_{k-1},a_\nu],
        \]
        
        \item[$(\beta)$] if 
        $F\in\tau_{T [\bar{b},M^{\eps+1}\rest\tau(T)]}\setminus\tau_{T}$ 
        is a $k$--place function symbol, then there is a $(k+1)$--place 
        function symbol $F_*\in\tau(\Phi_{\eps+1})\setminus\tau(\Phi_\eps)$ 
        such that
        \[
        \hspace{-1.2cm}
        M^{\eps+1}_{\bar b}\models\mbox{``}F[c_0,\ldots, c_{k-1}] = c
        \mbox{''}\quad\mbox{iff}\quad M^{\eps+1}\models \mbox{``}F_*[c_0,
        \ldots,c_{k-1},a_\nu]=c\mbox{''}.\]
    \end{enumerate}
\end{enumerate}

Let us carry out the induction; note that there is a redundancy
in our contraction: each relevant $\bar b$ is taken care of
in the $\eps$-th stage for every
$\eps<\kappa$ large enough, independently,
for the different
$\eps$-s.
\medskip

\noindent\underline{For $\eps=0$}:\\
Let $\Phi_0=\Phi$.
\medskip

\noindent\underline{For a limit $\eps$:}\\
Let $\Phi_\eps$ be the direct limit of $\LL\Phi_\zeta:\zeta<
\eps\RR$.
\medskip

\noindent\underline{For $\eps=\zeta+1$}:\\
Let the family of sequences of the form $\bar{\sigma}=\LL\sigma_i(x):i<
i^*\RR$, where $\sigma_i(x)$ is a unary term in $\tau(\Phi_\zeta)$,
$i^*<\kappa$, be listed as $\LL\bar{\sigma}^\gamma(x):\gamma<\theta
\RR$, with $\bar{\sigma}^\gamma(x)=\LL\sigma^\gamma_i(x):i<i_\gamma
\RR$. Let $M^*_\eps$ be a $\theta^+$--resplendent (hence strongly
$\theta^+$-homogeneous and $\kappa$-resplendent) elementary extension of
${\rm GEM}_{\tau (T)}({}^{\kappa\geq}\mu,\Phi_\zeta)$, and let $M_\eps=
M^*_\eps\rest\tau_{T}$, and choose $\nu^*\in{}^\kappa\mu$.
For each $\gamma<\theta$ let $\bar{b}^\gamma_{\nu^*}=:\LL\sigma^\gamma_i
(a_{\nu^*}):i<i_\gamma\RR$. Now, $(M_\eps,\bar{b}^\gamma_{\nu^*})$ can
be expanded to a model $M^\zeta_{\bar{b}^\gamma_{\nu^*}}$ of
$T^*[\bar{b}^\gamma_{\nu^*},M_\eps]$, and let
\[\tau(T^*[\bar{b},M_\eps])\setminus\tau_{T}=\{
R^{\eps,\gamma}_{j,n}:j<\theta, n<\omega\}\cup\{F^{\eps,
\gamma}_{j,n}:j<\theta,n<\omega\},\]
where $R^{\eps,\gamma}_{j,n}$ is an $n$--place predicate and
$F^{\eps,\gamma}_{j,n}$ is an $n$--place function symbol. Next we
shall define an expansion $M^+_\eps$ of $M^*_\eps$. Its
vocabulary is
\[
\tau(\Phi_\zeta) \cup \{R_{\eps,\gamma,j,n},F_{\eps,\gamma,j,n} :
j < \theta, n < \omega\},
\]
where $R_{\eps,\gamma,j,n}$ is an $(n+1)$--place predicate,
$F_{\eps,\gamma,j,n}$ is an $(n+1)$--place function symbol, and no
one of them is in $\tau(\Phi_\zeta)$ (and there are no repetitions in their
list).

Almost lastly, for $\nu\in {}^\kappa\mu$ let $g_\nu$ be an automorphism  of
$M_\eps$ mapping $\mathrm{GEM}_{\tau (T)}(\{\nu^*\},\Phi_\zeta)$ onto $\mathrm{GEM}_{\tau (T)}(\{\nu\},\Phi_\zeta)$; 
moreover such that for any $\tau(\Phi_\zeta)$--term $\sigma(x)$ we have
$g_\nu(\sigma(a_{\nu^*}))=\sigma(a_\nu)$ (hence 
$\xi < \kappa\ \Rightarrow\ g_\nu(a_{\nu^* \rest \xi}) = a_{\nu\rest\xi}$ 
using $\sigma(x)=F^*_\xi(x)$).

Now we actually define $M^+_\eps$ expanding $M^*_\eps$:
\[\begin{array}{r}
R^{M^+_\eps}_{\eps,\gamma,j,n}=\big\{(g_\nu(c_0),g_\nu(c_1),
\ldots,g_\nu(c_{n-1}),g_\nu(a_{\nu^*})):\quad\\
M^\zeta_{\bar{b}^\gamma_{\nu^*}}\models R^{\eps,\gamma}_{j,n}(c_0,
\ldots,c_{n-1})\big\},
  \end{array}\]
$F^{M^+_\eps}_{\eps,\gamma,j,n}$ is an $(n+1)$--place function
such that
\[\begin{array}{l}
M^\zeta_{\bar{b}^\gamma_{\nu^*}}\models F^{\eps,\gamma}_{j,n}(c_0,
\ldots,c_{n-1})=c\quad\mbox{ implies}\\
\ \\
F^{M^+_\eps}_{\eps,\gamma,j,n}(g_\nu(c_0),\ldots,g_\nu(c_{n-
1}),a_\nu)=g_\nu(c).
  \end{array}\]
We further expand $M^+_\eps$ to $M^{++}_\eps$, with vocabulary
of cardinality $\leq\theta$ and with Skolem functions.

Now we apply ``${}^{\kappa\geq}\mu$ has the Ramsey property''
(see  \cite[1.14(4)]{Sh:E59} 
see ``even'' there, 
\cite[1.18]{Sh:E59}) to get $\Phi_\eps=\Phi_{\zeta+1}$,
$\tau(\Phi_\eps)=\tau(M^{++}_\eps)$, such that for every $n<
\omega$, $\nu_1,\ldots,\nu_n\in{}^\kappa\mu$, and first order formula
$\varphi(x_1,\ldots,x_n)\in\bbL(\tau(\Phi_\eps))$, for some
$\eta_1,\ldots,\eta_n\in {}^\kappa\mu$ we have
\begin{enumerate}
    \item[$(\alpha)$] $M^{++}_\eps\models\varphi[a_{\eta_1},\ldots,a_{\eta_n}]$ iff 
    ${\rm GEM}_{\tau (T)}({}^{\kappa\geq}\mu,\Phi_\eps)\models \varphi[a_{\nu_1},\ldots,a_{\nu_n}]$,
    
    \item[$(\beta)$] $\LL\eta_1,\ldots,\eta_n\RR,\LL\nu_1,\ldots,\nu_n\RR$ are
    similar in ${}^{\kappa\geq}\mu$ (see \cite[VII]{Sh:c} or \ref{z20}).
\end{enumerate}
It is easy to check that $\Phi_\eps = \Phi_{\zeta+1}$ is as required.
\medskip

So we have defined the sequence $\LL\Phi_\eps : \eps < \kappa \RR$ satisfying 
the requirements above, and let $\Phi_\kappa$ be its limit. It is as required in the claim.
\end{PROOF}

\begin{Claim}
\label{1.24x}
Assume $\alpha\leq\kappa$ is a limit ordinal and 
$\LL \gm_\gamma, \bfh^\beta_\gamma : \gamma < \beta < \alpha\RR$ 
is an inverse system of approximations.
\begin{enumerate}[(1)]
    \item There are $\gm_\alpha$, $\bfh^\alpha_\gamma$
    (for $\gamma<\alpha$) such that 
    $\LL \gm_\gamma, \bfh^\beta_\gamma : \gamma < \beta < \alpha + 1\RR$ 
    is an inverse system of approximations continuous at $\alpha$.
    
    \item For the following properties, if each $\gm_{\gamma+1}$ 
    (for $\gamma<\alpha$) satisfies the property, \underline{then} 
    so does $\gm_\alpha$: term closed, semi full, almost full, resplendent, 
    weakly/strongly indiscernible, weakly/strongly nice, $\bfE$-indiscernible 
    for any $\bfE \in \bbE$, weakly/strongly good, weakly/strongly excellent.
\end{enumerate}
\end{Claim}

\begin{PROOF}{\ref{1.24x}}
 Let $\bfW_{\gm_\alpha}=\bigcup\limits_{\beta<\alpha}\bfW_{\gm_\beta}$, 
 and let $M_\beta=M_{\gm_\beta}$ for $\beta<\alpha$. We shall define 
 $\cF_\alpha=\cF_{\gm_\alpha}$, $M_\alpha=M_{\gm_\alpha}$ and 
 $N^\alpha_W=N^{\gm_\alpha}_W$ and 
 $M^\alpha_{\bar{b}} = M^{\gm_\alpha}_{\bar{b}}$ below.
\medskip

First let $\cF_\alpha$ (formal set, consisting of function symbols not
of functions), $\bfh^\alpha_\beta$ ($\beta<\alpha$) be the inverse limit of
$\LL\cF_\beta, \bfh^\beta_\gamma:\gamma\leq\beta<\alpha\RR$, i.e.,
\begin{enumerate}
    \item[$(\alpha)$] $\bfh^\alpha_\beta$ is a partial function from $\cF_\alpha$ 
    onto $\cF_\beta$ as in Definition \ref{1.22y}.
    
    \item[$(\beta)$] $\bfh^\alpha_\gamma = \bfh^\beta_\gamma \circ \bfh^\alpha_\beta$ for
    $\gamma<\beta<\alpha$,
    
    \item[$(\gamma)$] $\cF_\alpha = \bigcup\limits_{\beta<\alpha}\Dom(\bfh^\alpha_\beta)$,
    
    \item[$(\delta)$] If $\beta_*<\alpha$, $f_\beta\in \cF_\beta$, for
    $\beta\in [\beta_*,\alpha)$, satisfy $\bfh^\beta_\gamma(f_\beta) = f_\gamma$ when
    $\beta_*\leq\gamma<\beta<\alpha$, \underline{then} for one and only one
    $f\in \cF_\alpha$ we have:
    \[
    \zeta_f=\zeta_{f_\beta}\mbox{ for }\beta\in [\beta_*,\alpha)\qquad
    \mbox{ and }\qquad \eta_{f,\zeta} = \bigcup\left\lbrace\eta_{f_\beta,\zeta}:
    \beta_* \leq \beta<\alpha\right\rbrace,
    \]
    
    \item[$(\eps)$] every $f\in \cF_\alpha$ has the form of $f$ in $(\delta)$,

    \item[$(\zeta)$] $f^\ast_{\rho,\zeta}$ are as in (B) of Definition
    \ref{1.21}, i.e., for any $\rho\in {}^{\alpha}\mu$ and $\zeta < \kappa$
    we have $\beta < \alpha\ \Rightarrow\ \bfh^\alpha_\beta (f^*_{\rho\rest\beta}, \zeta) = f^{*,\gm_\beta}_{\rho,\zeta}$.
\end{enumerate}
\medskip

Second, we similarly choose $\bff_{\gm_\alpha}$.
\medskip

Thirdly, we choose $M_\alpha$ and interpretation of $f$ (for $f\in \cF_\alpha$) and $M^+_{\bar{b}}$ when
\[\bar{b}\in\{\Rang(f): f\in \cF_\alpha\ \&\ (\forall\zeta<\zeta_f)
(\exists\nu\in W)(\eta^f_\zeta\lhd\nu)\}\]
for some $W\in \bfW_{<\alpha}$. Though we can use the compactness theorem, it
seems to me more transparent to use ultraproduct . So let $D$ be an ultrafilter on
$\alpha$ containing all co-bounded subsets of $\alpha$. Let
$M_\alpha=\prod\limits_{\beta<\alpha} M_\beta/D$. If $f\in \cF_\alpha$,
let $\beta_f<\alpha$ and $\LL f_\gamma:\gamma\in [\beta_f,\alpha)
\RR$ be such that $\beta_f\leq\gamma<\alpha\ \Rightarrow\
\bfh^\alpha_\gamma(f)=f_\gamma$, so $\LL\eta^f_\zeta\rest\beta_f:
\zeta<\zeta_f\RR$ has no repetitions. Now, when $\eta^f_\zeta
\lhd\nu_\zeta\in {}^\kappa\mu$, let
\[f_\gm(\ldots,\nu_\zeta,\ldots)=\LL c_\gamma:\gamma<\alpha\RR/D,\]
where
\[
\begin{array}{rcl}
\gamma\in (\beta_f,\alpha)&\Rightarrow& c_\gamma=
\big(\bfh^\alpha_\gamma(f)_{\gm_\gamma}
\big)(\ldots,\nu_\zeta,\ldots)\in M_\gamma,\\
\gamma \leq \beta_f & \Rightarrow & c_\gamma \mbox{ is any member of } M_\gamma.
\end{array}
\]
So $M^\alpha_W$ is well defined for $W\in\bfW_{\gm(\alpha)}$.
\medskip

Fourth, if $\bar{b}=\LL b_\eps:\eps<\eps(*)\RR
\in{}^{\kappa>}(M^\alpha_W)$, $b_\eps=
f^\gm_\eps(\ldots,\nu_\zeta,
\ldots)_{\zeta<\zeta_{f_\eps}}$, and $\beta_*<\alpha$ and for
$\gamma\in [\beta_*,\alpha)$: $f_{\gamma,\eps}\in \cF_{\gm_\beta}$, $\LL f_{\gamma,\eps}:\eps<\eps(*)\RR\in\Dom({\bf
f}_{\gm_\beta})$, and $\bfh^\alpha_\gamma(f_{\gamma,\eps})=f_\eps$,
\underline{then} we let $\bar b^\beta=\LL b^\beta_\eps:
\eps<\eps(*)\RR$ where $b^\beta_\eps$ is
$f^{\gm_\beta} (\ldots,\nu_\zeta,\ldots)_{\zeta<\zeta_{f_\eps}}$
if $\beta\in [\beta_*,\alpha)$ and $b^\beta_\eps$ is any member of
$M_\gamma$ if $\beta<\beta_*$ and lastly
we define $M^\alpha_{\bar{b}}=\prod\limits_{\beta\in
[\beta_*,\alpha)} M^\beta_{\bar{b}_\beta}/D$.  We still have to check that
if for the same $\bar{b}$ we get two such definitions, then they agree, but
this is straightforward.
\medskip

Fifth, we choose $M^\alpha_{\bar{b}}$ for other $\bar{b}\in {}^{\kappa>}(
M_\alpha)$ for which $M^\alpha_{\bar{b}}$ is not yet defined
to satisfy clause (d) of Definition \ref{1.21}; note that by the
choice of $\bfW_{\gm_\alpha}$ those choices do not
influence the preservation of weakly/strongly indiscernible. So
$\gm_\alpha$ is well defined and one can easily check that it is as
required.
\end{PROOF}

\begin{Claim}
\label{1.25}
Assume $\alpha = \beta + 1 < \kappa$, and $\gm_1$ is a
$\beta$--approximation.
\begin{enumerate}
    \item There are $\bfh_*$ and an $\alpha$--approximation $\gm_2$ 
    such that $\gm_1 \leq_{\bfh_*} \gm_2$,\\ $M_{\gm_1}\prec M_{\gm_2}$,
    $M^{\gm_2}_{\bar{b}}=M^{\gm_1}_{\bar{b}}$, and $\Dom(\bfh_*)=\cF_{\gm_2}$.
    
    \item If $\gm_1$ is weakly/strongly nice, then $\gm_2$ is weakly/strongly nice.
    
    \item If $\gm_1$ is weakly/strongly indiscernible , \underline{then} 
    $\gm_2$ is weakly/strongly indiscernible; simply for $\bfE$-indiscernible,
    $\bfE \in \bbE_\alpha$.
\end{enumerate}
\end{Claim}

\begin{PROOF}{\ref{1.25}}  (1) Should be clear.

Let $\alpha(\gm_2)=\alpha$, $\bfW_{\gm_2}=\bfW_{\gm_1}$, $M_{\gm_1}\prec M_{\gm_2}$ 
and $M^{\gm_2}_{\bar{b}}=M^{\gm_1}_{\bar{b}}$ for $\bar{b}\in {}^{\kappa>}(M_{\gm_1})$. 
Then let
$$
\begin{array}{ll}
\cF_{\gm_2}=\lbrace g_{f,h}:& f\in\cF_\beta, h\mbox{ is a function with domain } 
\{\eta_{f,\zeta} : \zeta < \zeta_f\}\\
&\mbox{satisfying }h(\eta_{f,\zeta})\in\mathrm{Suc}(\eta_{f,\zeta})=\{\eta_{f,\zeta}
\conc\LL\gamma\RR:\gamma<\mu\}\rbrace,
\end{array}
$$
where for $g=g_{f,h}$ we let $\zeta_g=\zeta_f$ and $\eta_{g,\zeta}=h(
\eta_{f,\zeta})$, and if $\nu_\zeta\in I_{\eta_{g,\zeta}}$ for $\zeta<
\zeta_g$ ($=\zeta_f$), then
\[g^{\gm_2}_{f,h}(\ldots,\nu_\zeta,\ldots)=f^{\gm_1}
(\ldots,\nu_\zeta,\ldots)\in M_{\gm_1}\prec M_{\gm_2}.\]  We
define $\bfh_*$ by:
\[\Dom(\bfh_*)=\cF_{\gm_2}\quad\mbox{and}
\quad \bfh_*(g_{f,h})=f.\]
Lastly let
\[\begin{array}{ll}
\Dom({\bff}_{\gm_2})=\big\{\LL g_\eps:\eps<\eps(*)\RR:&
\mbox{for some }\bar{f}=\LL f_\eps:\eps<\eps(*)\RR\in \Dom({\bf
f}_{\gm_1})\\
&\mbox{and a function } h \mbox{ with domain }\\
&\{\eta_{f_\eps,\zeta}:\zeta<
\zeta_{\bar{f}}\} \mbox{ i.e., does not depend on } \eps\\
&\mbox{ we have } \eps<\eps(*)\ \Rightarrow\ g_\eps=g_{f_\eps,h}\ \big\},
  \end{array}\]
and if $\bfh,\bar{f},\bar{g}=\LL g_{f_\eps,h}:\eps<\zeta_{\bar{f}} \RR
\in \Dom(\bff_{\gm_2})$ are as above, $\sigma(\bar{x})$ is a
$\tau[T,\eps(*)]$--term,
$\bar{x}= \LL x_\xi:\xi\in u\RR$, and $u$ is a finite subset of
$\eps(*)$
and ($\bff_{\gm_1} (\bar{f})) (\sigma(\bar{x}))=f$, then
($\bff_{\gm_2} (\bar{g})) (\sigma(\bar{x}))=g_{f,h}$.

Now check.
\medskip

\noindent 2), 3)\quad Easy. 
\end{PROOF} 

\begin{Definition}
\label{1.26}
\begin{enumerate}
\item For approximations $\gm_1,\gm_2$, let $\gm_1\leq^*
\gm_2$ mean that $\alpha(\gm_1)=\alpha(\gm_2)$ and $\gm_1\leq_\bfh \gm_2$ 
with $\bfh$ being the identity on $\cF_{\gm_1}\subseteq \cF_{\gm_2}$, and $\bfW_{\gm_1}\subseteq {\bfW}_{\gm_2}$ and $\bff_{\gm_1}\subseteq \bff_{\gm_2}$ 
(the last condition means that if $\bar f \in \Dom(\bff_{\gm_1})$ then
$\bar f \in \Dom (\bff_{\gm_2})$ and the function $\bff_{\gm_2}
(\bar f)$ is equal to the function $\bff_{\gm_1} (\bar f)$.
\item Let $\gm_1<^*\gm_2$ mean that
\begin{enumerate}
\item[(a)] $\gm_1\leq^*\gm_2$,
\item[(b)] if $\bar{f}\in \bbF_{\gm_1}$ then $\bar{f}\in\Dom({\bf
f}_{\gm_2})$.
\end{enumerate}
\end{enumerate}
\end{Definition}

\begin{Observation}
\label{1.27}
\begin{enumerate}
\item $\leq^*$ is a partial order, $\gm_1\leq^*\gm_1$, and
\[\begin{array}{l}
\gm_1<^* \gm_2\ \Rightarrow\ \gm_1\leq \gm_2,\quad\mbox{ and}\\
\gm_1\leq^*\gm_2<^*\gm_3\ \Rightarrow\ \gm_1<^*\gm_3,\quad\mbox{ and}\\
\gm_1<^*\gm_2\leq\gm_3\ \Rightarrow\ \gm_1<^*\gm_3.
  \end{array}\]
\item Each $\leq^*$--increasing chain of length $<\theta^+$ has a lub
(essentially its union). If all members of the chain are weakly/strongly
indiscernible, then so is the lub.
\item If $\LL\gm_\eps:\eps<\kappa\RR$ is
$<^*$--increasing \underline{then} its lub $\gm$ is resplendent and
$\eps<\kappa\ \Rightarrow\ \gm_\eps<^*\gm$.
So if each $\gm_\eps$ is weakly/strongly good then $\gm$ is
weakly/strongly excellent.
\end{enumerate}
\end{Observation}

\begin{PROOF}{\ref{1.27}}
Easy. \\
\end{PROOF} 
As a warm up.  

\begin{Claim}
\label{2.15A}
\begin{enumerate}
\item For any $\alpha$--approximation $\gm_0$ there is a full, term
closed $\alpha$--approximation $\gm_1$ such that $\gm_0\leq^*
\gm_1$.
\item If $\gm_0$ is an $\alpha$--approximation, \underline{then}
there is a $\alpha$--approximation $\gm_1$ such that $\gm_0<^*\gm_1$ and $\Dom(\bff_{\gm_1}) = \bbF_{\gm_0}$.
\end{enumerate}
\end{Claim}

\begin{PROOF}{\ref{2.15A}}
1)\quad Let $M_{\gm_1}=M_{\gm_0}$, and $M^{\gm_1}_{\bar{b}}=M^{\gm_0}_{\bar{b}}$ for $\bar{b}\in {}^{\kappa>}
(M_{\gm_0})$. Let $\bfW_{\gm_1}=\bfW_\alpha$, and let
$\LL\bar{\nu}_\gamma:\gamma<\gamma^*\RR$ list the sequences
$\bar{\nu}\in {}^{\kappa>}({}^\kappa\mu)$ such that $\LL
\nu_\zeta\rest\alpha:\zeta<\lh(\bar{\nu})\RR$ is without
repetitions and $\{\nu_\zeta:\zeta<\lh(\bar{\nu})\}\notin\bfW_{\gm_0}$. Let $\bar{\nu}_\gamma=\LL\nu_{\gamma,\zeta}:\zeta<\zeta^*_\gamma
\RR$ and define $\bar{\rho}_\gamma=:
\LL \nu_ {\gamma,\zeta} \rest \alpha:\zeta<
\lh(\bar{\nu}_\gamma)\RR$,
and $W_\gamma=:\{\nu_{\gamma,\zeta}:\zeta<
\zeta^*_\gamma\}$ for $\gamma<\gamma^*$. Let $\beta_\gamma=
{\rm otp} \{\gamma_1<\gamma: (\forall \gamma_2<\gamma_1)
(\bar \rho_{\gamma_2}\neq \bar \rho_{\gamma_1})\}$.

For each $W\in\bfW_\alpha\setminus \bfW_{\gm_0}$, let $M^{\gm_1}_W$ be an elementary submodel of $M_{\gm_1}$ of cardinality
$\theta$ such that
\[\begin{array}{rcl}
W^*_1\subseteq W\ \wedge\ W^*_1\in \bfW_{\gm_0}& \Rightarrow&
M^{\gm_1}_{W^*_1}\prec M^{\gm_1}_W\qquad\quad\mbox{ and}\\
\bar{b}\in {}^{\kappa>}\big(M^{\gm_1}_W\big)& \Rightarrow& M^{\gm_0}_{\bar{b}}\rest |M^{\gm_1}_W|\prec M^{\gm_0}_{\bar{b}}.
\end{array}\]
Let $\LL a_{W,i}:i<\theta\RR$ list the elements of $M^{\gm_1}_W$.
For $\beta<\beta_{\gamma^*}$ and $i<\theta$ we choose $f_{\beta,i}$ such that 
$$\gamma<\gamma^* \wedge \beta_\gamma=\beta\ \Rightarrow\ \zeta_{f_{\beta,i}}
= \lh(\bar \nu_\gamma) = \lh(\bar \rho_\gamma) = \lh(\bar \nu_{\beta_\gamma}),\ \eta(f_{\beta,i},\zeta)=\rho_{\gamma,\zeta}$$ 
and we define
$f^{\gm_1}_{\beta,i}$ by: if $\nu_\zeta\in I_{\rho_{\gamma,\zeta}}$ for $\zeta<\zeta_{f_{\beta,i}}$, and 
$\LL \nu_\zeta : \zeta < \zeta_{f_{\beta,i}}^*\RR = \bar{\nu}_\gamma$ then
$f^{m_1}_{\beta,i}(\ldots, \nu_{\gamma,\zeta},\ldots) = a_{W_\gamma,i}$.

\noindent Next, $\cF_{\gm_1}$ is almost 
$\cF_{\gm_0} \cup \{f_{\beta,i}:\beta<\beta_{\gamma^*},\ i<\theta\}$: 
we just have to term-close it. Lastly $\bff_{\gm_1}$ is defined as 
$\bff_{\gm_0}$ recalling that $\Dom(\bff_{\gm_1})$ is required just to be
a subset of $\bbF_{\gm_0}$.
\medskip

\noindent 2)\quad Also easy.
\end{PROOF}

Let $M^*$ be a $\|M_{\gm_0}\|^+$--resplendent elementary
extension of $M_{\gm_0}$. We define an $\alpha$--approximation $\gm_1$ as follows:
\begin{enumerate}
    \item[(a)] $\alpha_{\gm_1}=\alpha_{\gm_0}$, 
    $\bfW_{\gm_1} = \bfW_{\gm_0}$, $M_{\gm_1}=M^*$,
    
    \item[(b)] if $\bar{b}\in {}^{\kappa>}(M_{\gm_0})$, then 
    $M^{\gm_1}_{\bar{b}}$ is an elementary extension of $M^{\gm_0}_{\bar{b}}$,
    
    \item[(c)] $\bff_{\gm_1}\supseteq\bff_{\gm_0}$ and
    $\Dom(\bff_{\gm_1})=\bbF_{\gm_0}$,
    
    \item[(d)] $\cF_{\gm_1}=\cF_{\gm_0}$
    
    \item[(e)] if $(\bff_{\gm_1}(\bar{f}))(\sigma_\xi(\bar{x}^\xi)) = f$, $\eta(\bar{f},\zeta)\lhd\nu_\eta\in {}^\kappa\mu$, and
    \[
        \bar{b}=\LL f^{\gm_1}_\eps(\ldots,\nu_\zeta,\ldots)_{
        \zeta<\zeta_{\bar{f}}}:\eps<\eps_{\bar{f}}\RR\quad\mbox{
        and }\quad\bar{x}^\xi=\LL x^\xi_i:i\in u\RR,
    \]
    \underline{then}
    \[
        f^{\gm_1}(\ldots,\nu_\zeta,\ldots)_{\zeta<\zeta_{\bar{f}}}=
        \sigma^{M^{\gm_1}_{\bar{b}}}(\LL f_i(\ldots,\nu_\zeta,
        \ldots)_{\zeta<\zeta_{\bar{f}}}: i\in u\RR).
    \]
\end{enumerate}

\begin{Main Claim}
\label{1.28}
Assume $\gm_0$ is a weakly nice approximation. \underline{then}
there is a weakly good
approximation $\gm_1$ such that $\gm_0<^*\gm_1$ with
$\bfW_{\gm_1}=\bfW_{\gm_0}$.
\end{Main Claim}

\begin{PROOF}{\ref{1.28}}
By \ref{2.15A}(1)+(2) there is a full term closed $\gm_1$
such that $\gm_0<^*\gm_1$ and $\Dom(\bff_{\gm_1})=\bbF_{\gm_0}$.
We would like to ``correct''$\gm_1$ so that it is weakly indiscernible.
Let $\gm_2$ be an $\alpha_{\gm_1}$--approximation as guaranteed
in Claim \ref{2.16A} below, so it is good and reflecting we 
clearly see that $\gm_0\leq^*\gm_2$ and even $\gm_0<^*\gm_2$.
\end{PROOF} 

\begin{Main SubClaim}
\label{2.16A}
(1) Assume $\gm_0$ is a weakly nice $\alpha$--approximation and
$\gm_0 <^* \gm_1$ and $\Dom(\bff_{\gm_1})=\bbF_{\gm_0}$ and $\bfW_{\gm_1}$ 
is an ideal (that is closed under finite union).
\underline{then} there is a good $\alpha$--approximation $\gm_2$ such that:
\begin{enumerate}
\item[(a)] $\alpha_{\gm_2}=\alpha_{\gm_1}$, $\cF_{\gm_2}=\cF_{\gm_1}$, $\bff_{\gm_2}=\bff_{\gm_1}$, and $\bfW_{\gm_2}=\bfW_{\gm_1}$.
\item[(b)] $\gm_0 <^* \gm_2$;
\end{enumerate} 
We may add

\noindent (c) Assume:

\begin{enumerate}
    \item[$(\alpha)$] $n<\omega$ and $f_\ell\in \cF_{\gm_1}$,
    $\nu^\ell_\zeta\in I_{\eta(f_\ell,\zeta)}$ for $\zeta<\zeta_{f_\ell}$,
    $\ell<n$, and $\Delta$ is a finite set of formulas in $\bbL (\tau_{T})$
    
    \item[$(\beta)$]  $m<\omega$ and for $k<m$ we have $\bar{f}^k=\LL
    f^k_\eps:\eps<\eps_k\RR\in \Dom(\bff_{\gm_1})$ and $n_k<
    \omega$ and $g_{k,\ell}\in \cF_{\gm_2}$ (for $\ell<n_k$) satisfying
    \[
    \LL\eta(g_{k,\ell},\zeta): 
    \zeta<\zeta_{g_{k,\ell}}\RR=\LL\eta(
    \bar{f}^k,\zeta):\zeta<\zeta_{\bar{f}^k}\RR,
    \]
    and $\nu^\ell_{k,\zeta}\in I_{\eta(\bar{f}^k,\zeta)}$ for $\ell<n_k$,
    $\zeta < \zeta_{\bar{f}^k}$, and $\Delta_k$ is a finite set of formulas in
    $\bbL(\tau[\zeta_{\bar{f}},\tau(T)])$.
\end{enumerate}

\underline{then} we can find $\rho^\ell_\zeta$ for $\ell<n_k$, $\zeta<\zeta_{f_\ell}$
and $\rho^\ell_{k,\zeta}$ for $\ell<n_k$, $\zeta<\zeta_{\bar{f}^k}$
for $\ell<n_k,k<m$ such that
\begin{enumerate}
    \item[(i)]   $\rho^\ell_\eps\in I_{\eta(f_\ell,\zeta)}$ for 
    $\zeta < \zeta_{f_\ell}$ and $\rho_{k,\zeta}\in I_{\eta(\bar{f}^k,\zeta)}$ for
    $\zeta < \zeta_{\bar{f}^k}$ for $\ell < n , k < m$,
    
\item[(ii)]  the sequences $\LL \rho^\ell_\zeta : \ell < n,\ 
\zeta < \zeta_{f_\ell}\RR \conc \LL \rho_{k,\zeta} : \ell < n_k , k < m,\ 
\zeta < \zeta_{\bar{f}^k} \RR$ and $\LL \nu^\ell_\zeta:\ell<n,\
\zeta<\zeta_{f_\ell}\RR\conc\LL \nu_{k,\zeta}:\ell<n_k,k<m,\ \zeta<
\zeta_{\bar{f}^k}\RR$ are similar\\ (see \cite[VII]{Sh:c} or \ref{z20}),

    \item[(iii)] the $\Delta$--type realized by the sequence
    \[
    \LL f^{\gm_2}_\ell(\ldots,\nu^\ell_\zeta,\ldots)_{\zeta<\zeta_{f_\ell}}:\ell<n\RR
    \]
    in $M_{\gm_2}$ is equal to the $\Delta$--type which the sequence
    \[
    \LL f^{\gm_1}_\ell(\ldots,\rho^\ell_\zeta,\ldots)_{\zeta < \zeta_{f_\ell}} : \ell < n\RR
    \]
    realizes in $M_{\gm_1}$,
    
    \item[(iv)]  for $k<m_1$, the $\Delta_k$--type realized by the sequence
    \[
    \LL g^{\gm_2}_{k,\ell}(\ldots,\nu^\ell_{k,\zeta},\ldots)_{\zeta < \zeta_{\bar{f}^k}} : \ell < n_k\RR
    \]
    in the model $M^{\gm_2}_{\LL f^{k,\gm_2}_\eps: \eps<\eps_k\RR} (\dots
    \nu^\ell_{k,\zeta}\ldots)_{\zeta<\zeta:\eps<\eps_k \RR}$
    is equal to the $\Delta_k$--type realized by the sequence
    \[\LL g^{\gm_1}_{k,\ell}(\ldots,\rho^\ell_{k,\zeta},\ldots)_{\zeta<
    \zeta_{\bar{f}^k}}:\ell<n_k\RR\]
    in the model $M^{\gm_1}_{\LL f^{k,\gm_2}_{\eps}:\eps<\eps_k\RR} (\dots
    \nu^\ell_{k,\zeta}\ldots)_{\zeta<\zeta:\eps<\eps_k\RR}$
    
    \item[(v)]   if $k_1,k_2<m$ then
    \[
    \begin{array}{l}
    \LL f^{k_1,\gm_2}_\eps(\ldots,\nu^{k_1}_\zeta,\ldots)_{\zeta<
    \zeta_{\bar{f}^{k_1}}}:\eps<\eps_{k_1}\RR=\\
    \LL f^{k_2,\gm_2}_\eps(\ldots,\nu^{k_2}_\zeta,\ldots)_{\zeta<
    \zeta_{\bar{f}^{k_2}}}:\eps<\eps_{k_2}\RR
    \end{array}
    \]
    { if and only if  }
    \[\begin{array}{l}
    \LL f^{k_1,\gm_1}_\eps(\ldots,\rho^{k_1}_\zeta,\ldots)_{\zeta<
    \zeta_{\bar{f}^{k_1}}}:\eps<\eps_{k_1}\RR=\\
    \LL f^{k_2,\gm_1}_\eps(\ldots,\rho^{k_2}_\zeta,\ldots)_{\zeta<
    \zeta_{\bar{f}^{k_2}}}:\eps<\eps_{k_2}\RR,
    \end{array}
    \]
    
    \item[(vi)]  if $\ell<n_k$, $k<m$, $\ell^*<n$, then
    \[
    f^{\gm_2}_\ell(\ldots,\nu^{\ell^*}_\zeta,\ldots)_{\zeta<\zeta_{f_\ell}}
    = f^{\gm_2}_{k,\ell}(\ldots,\nu^\ell_{k,\zeta},\ldots)_{\zeta<\zeta_{
    \bar{f}^k}}
    \]
    if and only if
    \[
    f^{\gm_1}_\ell(\ldots,\rho^{\ell^*}_\zeta,\ldots)_{\zeta<\zeta_{f_\ell}}
    = f^{\gm_1}_{k,\ell}(\ldots,\rho^\ell_{k,\zeta},\ldots)_{\zeta<\zeta_{
    \bar{f}^k}}.
    \]
\end{enumerate}
\end{Main SubClaim}

\begin{Discussion}
\label{1.29}
{\rm
Now we have to apply the Ramsey theorem to recapture weak
indiscernibility. Why do we only promise $\gm_0<^*\gm_1$ and 
$\Dom(\bff_{\gm_1}) = \bbF_{\gm_0}$, not that $\gm_1$ is excellent? 
Because $T^*[\bar{b},M]$ is not a continuous function of 
$(\bar{b},M)$, and more down to earth, during the proof we need to 
know the type of $\bar{b}$ whenever we consider types in
$M^{\gm_1}_{\bar{b}}$ in order to know $T^*[\bar b, M_\gm]$.
}
\end{Discussion}

Usually a partition theorem on what we already have is used at this moment,
but partition of infinitary functions tend to contradict ZFC. However, in
the set $\forGa$ expressing what we need, the formulas are finitary. So
using compactness we will reduce our problem to the consistency of the
set $\forGa$ of first order formulas in the variables
\[
    \big\{f(\ldots,\eta_\zeta,\ldots)_{\zeta<\zeta(f)} : f\in \cF^{\gm_1}
    \mbox{ and }\zeta < \zeta_f\ \Rightarrow\ \eta(f,\zeta) \lhd
    \eta_\zeta \in {}^{\kappa}\!\mu \big\}.
\]
This can be easily reduced to
the consistency of a set $\forGa$ of formulas in
$\bbL (\tau_{T})$
(first order).

We can get $\forGa$ because for all relevant $\bar b$ we know
$T^*[\bar{b},M]$.

\begin{PROOF}{\ref{2.16A}}
Let $Y=\{y_{f(\ldots,\nu_\eta,\ldots)_{\eta\in w[f]}}: f\in
\cF_{\gm_1}\ \mbox{ and }\ \nu_\eta\in I_\eta \mbox{ for } \eta\in w[f]\}$
be a set of
individual variables with no repetitions, recalling that
$w[f]=\{\eta [f,\eps]:
\eps<\zeta_f\}.$ For each $\bar{f}\in
\Dom(\bff_{\gm_1})$ and $\bar{\nu}=\LL\nu_\eta:\eta\in
w[\bar{f}] \RR$ such that $\nu_\eta\in I_\eta$, let
$\tau_{\bar{f},\bar{\nu}}$ be
$\tau[T,\lh \bar f]$ where $w[\bar{f}]=w[f_\eps]$ for each
$\eps<\lh(\bar{f})$; pedantically a copy of it over
$\tau_{T}$ so $(\bar{f}_1,\bar{\nu}_1)\neq (\bar{f}_2,\bar{\nu}_2)\
\Rightarrow\ \tau_{\bar{f}_1,\nu_1}\cap \tau_{\bar{f}_2,\bar{\nu}_2}=
\tau_{T}$. Let $\tau^*=\bigcup\{\tau_{\bar{f},\bar{\nu}}:
\bar{f},\bar{\nu}\mbox{ as above }\}\cup\tau_{T}$.

Let $\bfg_{\bar f, \bar \nu}$ be a one to one function from
$\tau[T, \lh (\bar f)]$ onto $\tau_{\bar f,\bar \nu}$ which is the identity
on $\tau_T$ preserve the arity and being a predicate function symbol,
individual constant. Let $\hat{\bfg}_{\bar f,\bar \nu}$ be the
mapping from $\bbL (\tau[T,\lh (\bar f)])$ onto
$\bbL(\tau_{\bar f,\bar \nu})$ which $\bfg_{\bar f,\bar \nu}$ induce.

We now define a set $\forGa$ (the explanations are for the use in the proof of
$\boxtimes_1$ below).

$\boxtimes_0$ \quad $\forGa=\forGa_0\cup\forGa_1\cup\forGa_2\cup\forGa_3\cup
\forGa_4\cup\forGa_5\cup \forGa_6 \cup \forGa_7 \cup \forGa_8$ where
\begin{enumerate}
\item[(a)]
$\forGa_0=\{\varphi_\zeta(y_{f^*_{\rho_1,\kappa}(\nu_1)},
y_{f^*_{\rho_2,\zeta+1}(\nu_2)})^\bft$: where $\bft=$ truth \underline{iff}\

$[\nu_1\rest (\zeta+1)=\nu_2\rest (\zeta+1)]$ and $\zeta<\kappa,\rho_\ell
\in {}^\alpha\mu$, and $\nu_\ell\in I^\alpha_{\rho_\ell}$
for $\ell=1,2\}$

 [explanation: to satisfy (iii) in clause (B) of Definition
\ref{1.21}].
\item[(b)]$\forGa_1=\{y_{f^*_{\rho_1,\zeta}(\nu_1)}=
y_{f^*_{\rho_2,\zeta}(\nu_2)} : \zeta < \kappa, \rho_\ell 
\in {}^\alpha\mu, \nu_\ell\in I^\alpha_{\rho_\ell}$
for $\ell=1,2$ and $\nu_1\rest \zeta = \nu_2\rest \zeta
\}$

[explanation:to satisfy (ii) in clause (B) of Definition \ref{1.21}].
\item[(c)] $\forGa_2=\{y_{f(\ldots\nu_\eta,\ldots)_{\eta\in w[f]}}=\sigma
(\ldots,y_{f_\ell(\ldots,\nu_\eta,\ldots)_{\eta\in w[f]}}\ldots):
f,\LL f_\ell:\ell<n\RR$ and $\LL\nu_\eta:\eta\in w[f]\RR$
are as in clause (E) of Definition \ref{1.22}(4) for $\gm_1\}$.

[explanation:this is preservation of the witnesses for closure under terms
of $\tau$,
in clause (E) of Definition \ref{1.22}(4) for $\gm_1$].
\item[(d)] $\forGa_3=\{y_{f(\ldots,\nu_\eta,\ldots)_{\eta\in w[f]}}=\sigma
(\ldots,f^\ell (\ldots,\nu_\eta (f^\ell,\zeta),\ldots)_{\zeta<\zeta(f^\ell)},
,\ldots)_{\ell<n} : f$, $\LL f^\ell:\ell<n\RR$ and 
$\LL f_\eps : \eps < \eps(*) \RR \in \Dom(\bff_{\gm_0})$ are as in clause (F) 
of Definition \ref{1.22}(4)$\}$.

[explanation: this is preservation of the witness for closure under terms
of the $\tau(M_{\bar b})$-'s
as in clause (F) of Definition \ref{1.22}(4) for $\gm_1$).
\item[(e)]$\forGa_4=\{\varphi(\ldots,y_{f_\ell(\ldots,\nu_\zeta,
\ldots)_{\zeta<\zeta^*}},\ldots)_{\ell<n}: \varphi(x_0,\ldots,x_{n-1}) \in
\bbL (\tau_{T})$ and $f_\ell\in \cF_{\gm_0}$ and
$\zeta^* = \zeta_{f_\ell}$ for $\ell<n$, and $\nu_\zeta\in
I_{\eta(f_\ell,\zeta)}$
for $\zeta<\zeta^*$ and $M_{\gm_0}\models\varphi[\ldots,
f_\ell^{\gm_0}(\ldots,\nu_\zeta,\ldots)_{\zeta<\zeta^*},
\ldots)_{\ell<n}\}$

[explanation: this is for being above $\gm_0$, the
$\bbL (\tau_{T})$-formulas].
\item[(f)] $\forGa_5$ like $\forGa_4$ for the $M_b$-'s that is

$ \forGa_5=\{\varphi(\ldots,y_{f_\ell(\ldots,\nu_\zeta,\ldots)_{\zeta<\zeta^*},
\ldots)_{\ell<n}}$: for some $\bar f,\bar \nu$ and $\LL f_\ell:\ell<n\RR$
we have $\varphi(x_0,\ldots, x_{n-1})\in \bbL (\tau_{\bar f,\bar \nu})$ and
$\bar{f} \in \Dom(\bff_{\gm_0}),
\zeta(\bar{f})=\zeta^*=\zeta_{f_\ell},f_\ell \in \cF_{\gm_0},
\eta(\bar{f},\zeta)=\eta(f_\ell,\zeta)$ for $\zeta<\zeta^*$, $\ell<n$
and $M_{\gm_0} \models \varphi [\ldots f^\gm_\ell
(\ldots,\nu_\zeta,\ldots)_{\zeta<\zeta^*} \ldots)_{\ell<n\}}$

[explanation: this is for being above $\gm_0$, the formulas from the
$M^{\gm_0}_{\bar{b}}-'s$].
\item[(g)] $\forGa_6=\{\varphi (\ldots,y_{f_\ell(\ldots\nu^1_\zeta,
\ldots)_{\zeta<\zeta^*}},\ldots)_{\ell<n} \equiv \varphi(\ldots,y_{f_\ell(
\ldots,\nu^2_\zeta,\ldots)_{\zeta<\zeta^*}},\ldots)_{\ell<n}:
\varphi (x_0,\ldots, x_{n-1})\in \bbL ))\tau_T)$ and $\zeta^*=\zeta(
f_\ell)$, $f_\ell\in \cF^{\gm_1},\nu^k_\zeta \in I_{\eta(f_\ell,
\zeta)}$ if $k=1,2$ and $\zeta<\zeta^*$ such that for exactly one
$\zeta<\zeta^*$ we have $\nu^1_\zeta \not= \nu^2_\zeta\}$

[explanation: this is for weak indiscernibility, 
see Definition \ref{1.23}(1)
clause (b) and Definition \ref{1.23}(3).].
\item[(h)] $\forGa_7=\{(\forall x_1\ldots x_n) [\varphi_1 (x_0,\ldots
x_{n-1})  
\equiv \varphi_2(x_0,\ldots,x_{n-1})]$:  
 for some $\bar f^1,\bar f^2
\in \bbF_\gm,\lh (\bar f^1)=\lh (\bar f^2)$, $\eta(\bar f^\ell,\zeta)
\triangleleft \nu^\ell_\zeta \in {}^\kappa\mu$ for $\zeta<\lh (\bar f^\ell)$
and $\eta(\bar f^1,\zeta_1)=\eta(\bar f^2,\zeta_2)\Rightarrow \nu^1_{\zeta_1}
=\nu^2_{\zeta_2}$; $\bar 
\varphi_\ell \in  \bbL   
(\tau_{{\bar f}^\ell,{\bar \nu}^\ell})$
and $\hat \bfg^{-1}_{\bar f^1,\bar \nu^1} (\varphi_1)=\hat \bfg^{-1}_{\bar f^2,\nu^2}
(\varphi_2)\}$.

[Explanation: this has to show the existence of the $M_{\bar b}$: we can avoid
this if we change the main definition such that instead $M_{\bar b}$
we have $M_{\bar f,\bar \nu}]$
\item[(i)] $\forGa_8=\{\varphi(\ldots, y_{f_\ell(\ldots,\nu^1_\zeta,\ldots)_{\zeta<\zeta^*}},
\ldots)_{\ell<n} \equiv \varphi (\ldots,
y_{f_\ell(\ldots,\nu^2_\zeta,\ldots)_{\zeta<\zeta^*}},\ldots)_{\ell<n}:$ for some
$\bar f\in \Dom (\bff_{\gm_1}), \bar f,\bar \nu,\bar \nu^1,\bar \nu^2$ and $\LL f_\ell:\ell<n\RR$ we have
$\varphi(x_0,\ldots,x_{n-1})\in \bbL (\tau_{\bar f,\bar \nu})$, $f_\ell \in
\cF_{\gm_1}, \zeta_{f_\ell}=\zeta^*,\eta(f_\ell,\zeta)=\eta_\zeta \in
{}^\alpha\mu$ for $\zeta<\zeta^*$, and $w(\bar f)\subseteq \{\eta_\zeta:\zeta<
\zeta^*\}$, and for $k=1,2$ we have $\eta_\zeta \triangleleft
\nu^k_\zeta \in {}^\kappa\mu$ for $\zeta<\zeta^*$,
and $\eta(\bar f,\zeta_1)=\eta_{\zeta_2} \Rightarrow \nu_{\zeta_1}=
\nu_{\zeta_2}^k\}$
\end{enumerate}
Clearly (e.g. for the indiscernibility we use term closure)
\begin{enumerate}
\item[$(\boxtimes_1)$] $\forGa$ is a set of first order formulas in
the free variables from $Y$ and the vocabulary $\tau^*$ such that an
$\alpha$--approximation $\gm$ satisfying (i) below is as required if and
only if clause (ii) below holds where
\begin{enumerate}
\item[(i)]  $\cF_\gm=\cF_{\gm_1}$, $\bff_\gm=\bff_{\gm_1}$, $\bfW_\gm=\bfW_{\gm_1}$,
\item[(ii)] interpreting $y_{f(\ldots,\nu_n,\ldots)_{\eta\in w[f]}}\in
Y$ as $f^\gm(\ldots,\nu_\eta,\ldots)_{\eta\in w[f]}$ and the predicates
and relation symbols in each $\tau_{\bar{f},\bar{\nu}}$ naturally,
$M_\gm$ is a model of $\forGa$ or more exactly not $M_\gm$ but
the common expansion of the $M^\gm_{\bar{b}}$-'s for $\bar{b} \in\{
\LL f^\gm_\eps (\ldots,\nu_\zeta,\ldots)_{\zeta<\zeta^*}:\eps
<\eps^* \RR: \bar{f}=\LL f_\eps:\eps<\eps^*
\RR \in \Dom(\bff_{\gm_0})$, and
$\nu_\zeta \in I_{\eta(\bar{f},\zeta)}\}$.
\end{enumerate}
\end{enumerate}
So it is enough to prove
\begin{enumerate}
\item[$(\boxtimes_2)$] $\forGa$ has a model.
\end{enumerate}
We use the compactness theorem, so let $\forGa^a\subseteq\forGa$ 
be finite. We say that
$\nu\in {}^\kappa\mu$ appears in $\forGa^a$, if for some variable
$y_{f(\ldots,\nu_\eta,\ldots)_{\eta\in w[f]}}$ appearing as a free variable in
some $\varphi\in\forGa^a$ we have  
$\nu\in\{\nu_\eta:\eta\in w[f]\}$
or some formula in $\forGa^a$ belongs to $\bbL (\tau_{\bar f,\bar \nu})
\setminus \bbL (\tau_T)$. We may
also say ``$\nu$ appears in $\varphi$'', and/or ``$f(\ldots,\nu_\eta,
\ldots)_{\eta\in w[f]}$ appears in $\forGa^a$'' (or in $\varphi$).

Let $n^*_0=|\forGa^a|$. Now, for each $\eta\in {}^\alpha\mu$ the set of
$\nu\in I^\alpha_\eta$ appearing in $\forGa^a$, which we call
$J^\alpha_\eta$,
is finite but on
$\bigcup\limits_{\eta\in{{}^\alpha\!\mu}}J^\alpha_\eta$ we
know only that its cardinality is $<\kappa$. Note that, moreover,
$n^*_1=:\max\{|
J^\alpha_\eta|:\eta\in{}^\alpha\mu\}$ is well defined $<\aleph_0$ as well as
$m^*_0=|\forGa^a\cap \forGa_0|$. For
each $\eta\in {}^\alpha\mu$ we can find a finite set $\bfu_\eta\subseteq\kappa$
such that:

\begin{enumerate}
\item[$(\otimes)$]
\begin{enumerate}
\item[\ (i)] if $\nu_1\neq\nu_2\in J^\alpha_\eta$, \underline{then} 
$\min\{\zeta : \nu_1(\zeta) \neq \nu_2(\zeta)\} \in \bfu_\eta $  
\item[(ii)] if $\varphi_\zeta[f^\ast_{\rho_1,\kappa}(\nu_1), 
f^\ast_{\rho_2, \zeta}(\nu_2)]^{\bft}$ from clause (B) appears in $\forGa^a\cap\forGa_0$,
then $\zeta, \zeta + 1 \in \bfu_\eta$
\item[(iii)] $\alpha\in \bfu_\eta$. 
\item[(iv)]  $|u_\eta|\leq   
 ( n^*_1)^2 
+ 2m^*_0+1$. 
\end{enumerate}
\end{enumerate}
Clearly $n^*_2=\max\{|\bfu_\eta|:\eta\in{}^\alpha\mu\}$ is well defined
($<\aleph_0$), so without loss of generality, $\eta\in {}^\alpha\mu\
\Rightarrow\ |\bfu_\eta|=n^*_2$.

Let $v\subseteq {}^\alpha\mu$ be finite, 
 in fact of size
$\leq |\forGa^a|=n^*_0$ such that: 
\begin{enumerate}
\item[(I)] if
$\varphi_\zeta[y_{f^\ast_{\rho_1,\kappa}(\nu_1)},
y_{f^\ast_{\rho_2,\zeta}(\nu_2)}]^{\bft}$ appears in 
$\forGa^a\cap \forGa_0$, so
$\ell \in \{1,2\} \Rightarrow \nu_\ell\in J^\alpha_{\rho_\ell}$, 
then $\rho_\ell \in v$ for $ {\ell} \in \{ 1,2 \} $,  
\end{enumerate}

Now, for all $\eta\in {}^\alpha\!\mu \setminus v$ we replace in $\forGa^a$ all
members of $J^\alpha_\eta$ by one $\nu_\eta\in I^\alpha_\eta$ and we
call what we get $\forGa^b$, i.e., we identify some variables. It suffices
to prove $\forGa^b$ is consistent. Now, by the choice of the set $v$ also
$\forGa^b$ is of the right kind, i.e., $\subseteq\forGa$.

[Why? We should check the formulas $\varphi$, in $\forGa^a \cap \forGa_i$ for each $i\leq 8$;
let it be replaced by $\varphi'\in \forGa^b$. If in $\varphi\in \forGa_0\cap \forGa^a$ by clause (ii)
of $\otimes$ this substitution has no affect on $\varphi$. If $\varphi\in
\forGa_1,$ either $\varphi'=\varphi$ or $\varphi'$ is trivially true. If
$\varphi\in \forGa_3$, clearly $\varphi'\in \forGa_3$. If $\varphi\in
\forGa_4$ then $\varphi'\in \forGa_4$ as $\gm_0$ is nice hence weakly
indiscernible, i.e. clause (b) of  
Definition \ref{1.23}(1) (and the demand
$f_\ell\in \cF_ { \gm   _ 0 } $). 
If $\varphi \in \forGa_5$, similarly using clause
(c) of Definition \ref{1.23}(1). Lastly if $\varphi\in \forGa_6$ we just note
that similarly is preserved and similarly for $\varphi\in
\forGa_7\cup \forGa_8$].

We then transform $\forGa^b$ to $\forGa^c$ by
replacing each $\varphi$ by $\varphi'$, gotten by replacing, for each
$\rho\in v$, every $\nu\in J^\alpha_\rho$ by $\nu^{[*]}\in {}^\kappa\mu$
where $\nu^{[*]} (\beta)=\nu (\beta)$ if $\beta\in \alpha\cup u_\rho$
and $\nu^{[*]} (\beta)=0$ otherwise.
It suffices to prove the consistency of $\forGa^c$. Now, the effect is
renaming variables and again $\forGa^c\subseteq\forGa$.
Let $\rho^*=\LL \rho^*_k:k<k^*\RR$
list the $\rho\in {}^\kappa\mu$
which appear in $\forGa^c$ such that $\rho\rest \alpha\in v$. Let $\eta_k=
\rho^*_k\rest \alpha$ so $\eta_k\in v$, and let
$$
\begin{array}{ll}
    \Upsilon = \big\{\bar \rho : & \bar \rho = 
    \LL \rho_k : k < k^*\RR,\ 
    \eta_k \lhd \rho_k \in {}^\kappa\!\mu, \\
    & (\forall \eps) \big[\alpha\leq \eps < \kappa \wedge 
    \eps \notin u_{\eta_k}\Rightarrow \rho_k(\eps) = 0 \big]\\
    & \text{and $\bar \rho$ is similar to $\bar \rho^*$}\big\}
\end{array}
$$
(i.e., for $k_1, k_2<k^*$ and $\eps < \kappa$ we have
$\rho_{k_1}(\eps)<\rho_{k_2}(\eps) \Rightarrow
\rho^*_{k_1}(\eps)<\rho^*_{k_2}(\eps)$.)

For each $\bar \rho\in\Upsilon$ we can try the following model 
as a candidate to be a model of $\forGa^c$.
It expands $M_{\gm_1}$, and if symbols from 
$\tau_{\bar f,\bar \nu}\setminus \tau_T$ appear they are 
interpreted as their $\bfg^{-1}_{\bar f,\bar \nu}$-images are 
interpreted in $M^{\gm_1}_{\LL f_\eps(\bar \nu) : \eps < \lh(\bar f)\RR}$. 
Lastly we assign to the variable 
$y_f(\ldots,\nu_\zeta,\ldots)_{\zeta<\zeta_f}$ appearing
in $\forGa^c$  the element 
$f_{\gm_1}(\ldots,\nu_\zeta,\ldots)_{\zeta<\zeta_f}$ of $M_{\gm_1}$.
Call this the $\bar \rho$-interpretation.
Considering the formulas in $\forGa^c \cap \forGa_i$ for
$i\in \{0,\ldots, 5,7\}$ they always holds. For the formulas in
$\forGa^c\cap \forGa_6, \forGa_8$ we can use a partition
theorem on trees with $|n^*_2|<\aleph_0$ levels (use 
     \cite[1.16]{Sh:E59}(4),
which is an overkill, but has the same spirit 
(or \cite[AP2.6,p.662]{Sh:c})).
\end{PROOF} 

\begin{Claim}
\label{1.33}
There is an increasing continuous inverse
system of approximations
\[\LL\gm_\gamma,\bfh^\alpha_\beta:\gamma\leq\kappa,\ \beta\leq\alpha
\leq\kappa\RR\]
such that each $\gm_\gamma$ is weakly excellent.
\end{Claim}
 
\begin{PROOF}{\ref{1.33}}
By induction on $\alpha\leq\kappa$ we choose $\gm_\alpha$ and
$\LL \bfh^\alpha_\beta : \beta < \alpha\RR$ with our inductive hypothesis being
\begin{enumerate}
    \item[$(*)$]
    \begin{enumerate}
        \item[(a)]  
        $\LL\gm_{\beta_1}, \bfh^\beta_\gamma : 
        \beta_1\leq\alpha, \gamma < \beta \leq \alpha\RR$ 
        is an inverse system of approximations,
        
        \item[(b)] $\gm_\beta$ is a weakly excellent $\beta$--approximation.
    \end{enumerate}
\end{enumerate}

\noindent\underline{For $\alpha=0$}:\\
A weakly excellent good $0$--approximation exists by \ref{1.23x}.
\medskip

\noindent\underline{For $\alpha$ limit}:\\
Clearly $\LL\gm_{\beta_1},
\bfh^\beta_\gamma:\beta_1<\alpha,\gamma<
\beta<\alpha\RR$ is an inverse system of good
weakly excellent approximations with
$\alpha(\gm_\beta) = \beta$. So by \ref{1.24x} we can find 
$\gm_\alpha,\bfh^\alpha_\beta$ ($\beta<\alpha$) as required.
\medskip

\noindent\underline{For $\alpha=\beta+1$}:\\
By \ref{1.25}(1+2) there is $\gm^*_{\alpha,0}$ a weakly nice
$\alpha$--approximation such that $\gm_\beta\leq^* \gm^*_{
\alpha,0}$. By \ref{1.28} there is a full term closed $\alpha$--approximation
$\gm^*_{\alpha,1}$ such that $\gm^*_{\alpha,0} \leq^*
\gm^*_{\alpha,1}$ and $\gm^*_{\alpha,1}$ is good.
We can choose by induction on $\eps\in [1,\kappa]$ good
$\alpha$--approximations $\gm_{\alpha,\eps}$,
$\leq^*$--increasing continuously, $\gm_{\alpha,\eps}<^*\gm_{\alpha,\eps+1}$.

For $\eps=1$, $\gm_{\alpha,\eps}$ is defined; for
$\eps$ limit use \ref{1.27}(2), for $\eps$ successor use
\ref{1.28}, and $\gm_\alpha=:\gm_{\alpha,\kappa}$ is good by
\ref{1.27}(3).  
\end{PROOF}

\begin{Claim}
\label{1.34}
Assume $\gm_\alpha, \bfh^\alpha_\gamma$ for $\alpha\leq \kappa,
\gamma<\alpha$  
as in
\ref{1.33} with $\mu=\lambda$ and $\lambda=\lambda^\kappa\geq\theta$
(e.g., $\lambda=\lambda^\kappa\geq 2^{|T|}$). \underline{then} there
are $>\lambda$ pairwise non-isomorphic $\kappa$--resplendent models of
$T$ of cardinality $\lambda$.
\end{Claim}

\begin{PROOF}{\ref{1.34}}

 Let $\gm=\gm_\kappa$
 and 
$I \subseteq {}^{\kappa\geq}\lambda$,
 $|I|=\lambda$ and for
simplicity
$\{\eta\in {}^{\kappa}\lambda:\eta(\eps)=0$ for every large
enough $\eps<\kappa\}\cup
{}^{\kappa>}\lambda\subseteq I$. Let $M_I$ be the submodel of
$M_\gm$ with universe
$$
\left\lbrace f\left(\ldots,\nu_{\eta(f,\zeta)},\ldots\right)_{\zeta<
\zeta_f}:f\in \cF_\gm
\mbox{ and }\eta(f,\zeta)\in I\cap {}^\kappa
\lambda\mbox{ for every }\zeta<\zeta_f\right\rbrace.
$$
Trivially, $\|M_I\|\leq\lambda^\kappa=\lambda$ and by clause (B) of
Definition \ref{1.21} clearly
by \ref{1.17}(1) it follows that
the sequence 
$\LL a_\eta:\eta\in {}^{\eps}\lambda\RR$
is with no repetitions for each $\eps<\lambda$
 hence by the indiscernibility 
the sequence 
$\LL a_\eta:\eta\in I\RR$ is with no
 repetition, so
 $\|M_I\|\geq |I|\geq \lambda$, so $\|M_I\|=
\lambda$.

Now, $M_I$ is a $\kappa$-resplendent model of $T$ as $\gm$ being
weakly excellent is full and resplendent.

For $\zeta<\kappa$, $\nu\in {}^\zeta\lambda$ let
$a_\nu=f_{\eta,\zeta}^{*,\gm_\kappa}(\eta)$ ($\in
M_I$) for any $\eta\in I^\zeta_\nu \cap I$.

The point is:
\begin{enumerate}
\item[$(\otimes)$] For $\eta\in {}^\kappa\lambda$, $\nu_\gamma\in {}^{\gamma
+1}\lambda$, $\nu_\gamma\rest\gamma=\eta\rest\gamma$,
$\nu_\gamma\neq\eta\rest(\gamma+1)$, we have:
\begin{enumerate}
\item[$\circledast$] the type $
\left\lbrace\varphi (x,a_{\eta\rest(\gamma+1)})\equiv\lnot\varphi
(x,a_{\nu_\gamma}):\gamma<\kappa\right\rbrace$

is realized in $M_I$ \underline{iff}\ $\eta\in I$.
\end{enumerate}
\end{enumerate}
[Why? The implication ``$\Leftarrow$'' holds by clause (B)(iii)
 of Definition
\ref{1.21}. For the other direction, if $c\in M_I$, then for some
$W\in\bfW_\kappa$, satisfying $W\subseteq I$,
we have $c\in N^\gm_W$, and as
$\eta\notin
 I$ and $|W|<\kappa$ clearly for some $\alpha<\kappa$ we have
\[\{\nu:\eta\rest\alpha\lhd\nu\in{}^\kappa\mu\}\cap W=
\varnothing.\]
Let $c=f_{\gm_\kappa}
(\ldots,\nu_\zeta,\ldots)_{\zeta<\zeta_f}$, where $f\in \cF_\kappa$, so $\nu_\zeta=\eta(f,\zeta)$. By the continuity of the system,
for some $\gamma\in (\alpha,\kappa)$ we have
$f\in\Dom(\bfh^\kappa_\gamma)$,
and it suffices to prove that
\[M_I\models\mbox{`` }\varphi[c,a_{\eta\rest (\gamma+1)}]\equiv
\varphi[c,a_{\nu_\gamma}]\mbox{ ''.}\]
By the definition of a system, $\gm$ is full. Choose $\nu\in
I_{\nu_\zeta}$; recalling $\gm_\gamma\leq_{\bfh^\kappa_\gamma}
\gm_\kappa$ it suffices to prove that
\[\begin{array}{r}
M_{\gm_\gamma}\models\mbox{`` }
\varphi[(\bfh^\kappa_\gamma(f))(\ldots,
\nu_\zeta,\ldots)_{\zeta<\zeta_f}, f^*_{\eta\rest\gamma,\gamma}
(a_\eta)]\equiv\\
 \varphi[(\bfh^\kappa_\gamma(f))
(\ldots,\nu_\zeta,\ldots)_{\zeta<\zeta_f},
f^*_{\nu\rest\gamma,\gamma}(a_\nu)]\mbox{ ''.}
  \end{array}\]
But $\gm_\gamma$ is weakly excellent, hence it is ${\bf
E}^0_\gamma$--indiscernible, and hence the requirement holds.]

Now use 
\cite[\S2]{Sh:309}  
      to get, among those models,
$>\lambda$ which are non-isomorphic; putting in the 
$\eta\in {}^{\kappa}{\lambda}$ which are eventually zero does not matter.
\end{PROOF}


 \newpage

\section{Strengthening}

\begin{Claim}
\label{3.1}
If there is strongly excellent $\kappa$-approximation $\gm$ and $\mu\geq
\lambda=\lambda^\kappa\geq 2^{|T|}$, \underline{then} $T$ has
$2^\kappa$ non-isomorphic $\kappa$--resplendent models of cardinality
$\lambda$.
\end{Claim}

\begin{PROOF}{\ref{3.1}}
This time use Theorem 
\cite[2.3]{Sh:331}. For any $I \subseteq {}^{\kappa\geq}\lambda$ which includes
${}^{\kappa>}\lambda$, let $M_I\prec M_{\gm_\kappa}$ be defined
as in the proof of \ref{1.34}. For $\eta\in {}^\kappa\lambda$ let
$a_\eta=f^*_{\eta,\kappa}(\eta)$, and for $\eta\in {}^{\kappa>}\lambda$ of
length
$\gamma+1$ let $\eta'=\eta\rest\gamma\conc\LL\eta(\gamma)+1\RR$,
and for any $\nu\in I_\eta$, $\nu'\in I_{\eta'}$ let $\bar{a}_\eta=\LL
f^*_{\nu,\gamma+1}(\nu),f^*_{\nu',\gamma+1}(\nu')\RR$; the choice of
$(\nu,\nu')$ is immaterial. Let
\[\varphi(\LL\bar{x}_\alpha:\alpha<\kappa\RR)=(\exists
y)\left(\bigwedge_{
\alpha<\kappa}(\varphi(y,x_{\alpha,0})\equiv\neg
\varphi(y,x_{\alpha,1})\right).\]
 Now we can choose $f_I:M_I
\longrightarrow \clM_{\lambda,\kappa}$ such that
\smallskip
\begin{enumerate}
\item[(a)] if $f_I(b)=\sigma(\LL t_i:i<i^*\RR)$ such that  $t_i\in I\cap
{}^\kappa\lambda$ with no repetitions and $\sigma \in
\tau[\clM_{\lambda,\tau}]$,

\underline{then} for some $W\in\bfW_{\gm}$
and $\gamma<\kappa$ such that
$\LL \eta\rest \gamma:\eta\in W\RR$ is with no repetition we
have $\{t_i:i<i^*\}=W$ and for some $f\in \cF_\gm$ with
$\zeta_f=\zeta_*$, and $\eta(f,\zeta)=t_\zeta$ for $\zeta<\zeta^f$ we have
$b=f_\gm(\ldots, t_\zeta,\ldots)_{\zeta<\zeta_f}$ and
\item[(b)] $f_1(b)=\eta\in I$ if $b=a_\eta$ (see above).
\end{enumerate}
\smallskip

\noindent The new point is that we have to prove the statement $(*)$ in
     \cite[2.3]{Sh:331}(c)($\beta$).

So assume that for $\ell=1,2$ and $\alpha<\kappa$: $\bar{b}^\ell_\alpha\in
{}^2(M_{I_\ell})$, $f_{I_\ell}(\bar{b}^\ell_\alpha)=\bar{\sigma}_\alpha^\ell
(\bar{t}_\alpha^\ell)$, and $\bar{t}^1_\alpha=\LL t^2_{\alpha,\eps}:\eps<
\eps_\alpha\RR$. Assume furthermore that $\bar{\sigma}^1_\alpha=
\bar{\sigma}^2_\alpha$, $\bar{t}^1_\alpha=\bar{t}^2_\alpha$ for
$\alpha<\kappa$ (call it then $\bar{\sigma}^\alpha (\bar{t}_\alpha)$ though
possibly $I_1 \not= I_2$),
and the truth value of each statement
\[(\exists\nu\in I_\ell\cap{}^\kappa\lambda)\Big(\bigwedge_{i<\kappa}\nu
\rest\eps_i= t^\ell_{\beta_i,\gamma_i}\rest\eps_i\Big)\]
does not depend on $\ell\in\{1,2\}$. Assume further that $M_{I_1}\models
\varphi(\ldots,\bar{b}^1_\gamma,\ldots)_{\gamma<\kappa}$, and we shall prove
that $M_{I_2}\models \varphi(\ldots,\bar{b}^2_\gamma,\ldots)_{\gamma<
\kappa}$; this suffices.

First note that, as $f_{I_1},f_{I_2}\subseteq f_{({}^{\kappa\geq}\lambda)}$,
necessarily $\bar{b}^1_\alpha=\bar{b}^2_\alpha$ (so call it
$\bar{b}_\alpha$). Now, $M_{I_1}\models\varphi(\ldots,\bar{b}^1_\gamma,
\ldots)_{\gamma<\kappa}$ means that for some $c_1\in M_{I_1}$ we have
\[M_{I_1}\models\bigwedge_{\gamma<\kappa}\varphi[c_1,b_{\gamma,0}]\equiv\neg
\varphi[c_1,b_{\gamma,1}],\]
and let $c_1=f_1(\ldots,\eta,\ldots)_{\eta\in w[f_1]}$. Let
\[\begin{array}{lcl}
J&=&\{\eta:\eta\trianglelefteq t_{\alpha,j}\mbox{ for some }\alpha<\kappa,\
j<\lh(\bar{t}_\alpha)\;\}\qquad\mbox{ and}\\
J^+_\ell&=&\{\eta: \eta\in I_\ell\mbox{ or }\lh(\eta)=\kappa \mbox{ and }
(\forall\alpha<\kappa)(\eta\rest\alpha\in J)\;\}.
\end{array}\]
By the assumption, $J$ is $\lhd$--closed, $J\subseteq I_1\cap
I_2$, moreover $J^+_1=J^+_2$. Let $\gamma<\kappa$ be minimal such that
$\eta\in w[f_1]\setminus J^+\ \Rightarrow\ \eta\rest \gamma\notin J$,
and the sequence $\LL\eta(f_1,\zeta)\rest\gamma:\zeta<\zeta_f
\RR$ is with no repetitions and $f_1\in\Dom(\bfh^\kappa_\gamma)$.

Now we can choose $\nu_\eps\in I_{\eta(f_1,\zeta)\rest\gamma}$
from $I_2$
such that $\eta(f_1,\zeta)
\in J^+\ \Rightarrow\ \nu_\eps=\eta(f_1,\zeta)$. Let $f_2\in \cF_{\gm_\kappa}$ be such that $h^\kappa_\gamma(f_2)=h^\kappa_\gamma(
f_1)$ and $\eta(f_2,\zeta)=\nu_\zeta$ for $\zeta<\zeta_{f_2}=\zeta_{f_2}$.
Easily, $c_2=f^{\gm_\kappa}(\ldots,\nu_\zeta,\ldots)_{\zeta<
\zeta_{f_1}}\in M_{I_2}$ witness that
\[M_{I_2}\models(\exists y)[\bigwedge_{\alpha<\kappa}\varphi_\alpha(x,
b_{\alpha,0})\equiv\neg \varphi(x,b_{\alpha,1})]\]
(recalling $M_{I_1},M_{I_2}\prec M_{({}^{\kappa\geq}\lambda)}$).

\end{PROOF}

Recall and add

\begin{Definition}
\label{3.2}
\begin{enumerate}
\item $\bfE^1_\alpha\in\bbE$ (see Definition \ref{1.22x}) is defined like
$\bfE^0_\alpha$ (see Definition \ref{1.22x}(3)) except that we omit
clause (iv) there.
\item For $\alpha<\kappa$ define $\bfE^2_\alpha\in \bbE$ as
the following equivalence relation on $\{\bar \nu:\bar \nu\in
{}^{\kappa>}({}^\kappa\mu),\bar \nu$ with no repetition$\}$

$\bar \nu^1 \bfE^2_\alpha \bar \nu^2$ \underline{iff}\
\begin{enumerate}
\item[(i)] $\bar \nu^1,\bar \nu^2\in {}^{\kappa>}({}^\kappa\mu)$ are with no
repetition.
\item[(ii)] $\bar \nu^1, \bar \nu^2$ have the same length, all it $\zeta^*$.
\item[(iii)] $\nu^1_\zeta\rest \alpha=\nu^2_\zeta\rest \alpha$ for $\zeta<\zeta^*$.
\item[(iv)] for every $\zeta\in {}^\alpha\mu$, the sets $u^\ell_\eta=\{\zeta<\zeta^*:
\eta\triangleleft \nu^\ell_\zeta\}$ are finite equal and
$\LL \nu^1_\zeta:\zeta\in u^1_\eta\RR$, $\LL \nu^2_\zeta:\zeta\in u^2_\eta
\RR$ are  similar.
\end{enumerate}
\end{enumerate}
\end{Definition}

\begin{Claim}
\label{3.3}
\begin{enumerate}
\item In \ref{1.33} we can demand that every $\gm_\gamma$ is ${\bf
E}^1_\gamma$--indiscernible i.e. get the strong version.
\item Moreover we can get even $\bfE^2_\alpha$-indiscernibility.
\end{enumerate}
\end{Claim}

\begin{PROOF}{\ref{3.3}}
(1)  Very similar to the proof of \ref{1.33}. In fact, we need to repeat
\S 2 with minor changes. One point is that defining ``good'' we use ${\bf
E}^1_\gamma$; the second is that we should not that this indiscernibility
demand is preserved in limits, this is \ref{1.24x}, \ref{1.27}.
In fact this is the ``strongly'' version which is carried in \S2 the until
\ref{1.28}. From then on we should replace ``weakly'' by
``strongly'' and change the definition of $\forGa_6, \forGa_8$
appropriately in the proof of \ref{2.16A}.

(2)  Similarly, only we need a stronger partition theorem in the end of the
proof of \ref{2.16A}, but it is there anyhow.
\end{PROOF}

\begin{Remark}
{\rm
Clearly in many cases in \ref{3.1}, $\lambda=\lambda^{<\kappa}\geq\theta$
suffices, and it seems to me that with high probability for all.
Similarly  
for getting many $\kappa$--resplendent
models no one elementarily embeddable into another.
}
\end{Remark}

\newpage


\bibliographystyle{amsalpha}
\bibliography{shlhetal}

\end{document}